\documentclass[12pt]{article}
\usepackage{amsmath,amsfonts,amssymb,amsthm}
\title{\bf The rank two lattice type vertex operator algebras $V_L^+$ and their
automorphism groups}

\footnotetext{The first author is supported by NSF grants,
a NSF grant of China and a
research grant from UC Santa Cruz.}\footnotetext{
The second author is supported by NSA grant USDOD-MDA904-03-1-0098.    }
\author{Chongying Dong and  Robert L.~Griess Jr.}
\date{}

\makeatletter

\begin{document}

\newtheorem{thm}{Theorem}[section]
\newtheorem{prop}[thm]{Proposition}
\newtheorem{lem}[thm]{Lemma}
\newtheorem{rem}[thm]{Remark}
\newtheorem{coro}[thm]{Corollary}
\newtheorem{conj}[thm]{Conjecture}
\newtheorem{de}[thm]{Definition}
\newtheorem{notation}[thm]{Notation}
\newtheorem{nota}[thm]{Notation}
\newtheorem{hyp}[thm]{Hypothesis}
\newcommand{\fusion}[3]{{\binom{#3}{#1\;#2}}}

\def\labtt#1{\label{#1}}
\def\labttr#1{\label{#1}\rm }

\pagestyle{plain}
\maketitle

\def\vv{ {\hbox {\bf 1} }}  

\def \CC{\mathbb C}
\def \FF{\mathbb F}
\def \II{\mathbb I} 
\def \NN{\mathbb N}
\def \OO{\mathbb O}
\def \PP{\mathbb P}
\def \QQ{\mathbb Q}
\def \RR{\mathbb R}
\def\TT{\mathbb T}  
\def\UU{\mathbb U} 
\def \XX{\mathbb X}
\def \ZZ{\mathbb Z}

\def \D{{\cal D}}
\def \wt{{\rm wt}}
\def \tr{{\rm tr}}
\def \sp{{\rm span}}
\def \Res{{\rm Res}}
\def \Res{{\rm QRes}}
\def \End{{\rm End}}
\def \E{{\rm End}}
\def \Ind {{\rm Ind}}
\def \Irr {{\rm Irr}}
\def \Aut{{\rm Aut}}
\def \Hom{{\rm Hom}}
\def \mod{{\rm mod}}
\def \ann{{\rm Ann}}
\def \<{\langle} 
\def \>{\rangle} 
\def \t{\tau }
\def \a{\alpha }
\def \e{\epsilon }
\def \l{\lambda }
\def \L{\Lambda }
\def \g{\frak g}
\def \b{\beta }
\def \om{\omega }
\def \o{\omega }
\def \k{\kappa}
\def \c{\chi}
\def \ch{\chi}
\def \cg{\chi_g}
\def \ag{\alpha_g}
\def \ah{\alpha_h}
\def \ph{\psi_h}
\def \pf{{\bf Proof. }}
\def \voa{{vertex operator algebra\ }}
\def \svoa{{super vertex operator algebra\ }}
\def \qed{\mbox{$\square$} \vskip .6cm }
\def \lc{L_C}
\def \tlc{\widetilde{L}_C}
\def \tv{\widetilde{V}_L}
\def \vlc{V_{L_C}}
\def\tvlc{\widetilde{V}_{L_C}}
\def\vtlc{V_{\widetilde{L}_C}}
\def\tvtlc{\widetilde{V}_{\widetilde{L}_C}}
\def\ha{\frac{1}{2}}
\def\se{\frac{1}{16}}

\def\Ve{V^{0}}
\def\M{\rm \ {\sl  M}\llap{{\sl I\kern.80em}}\ }
\def\xx{\em}
\def\s{\sigma}
\def\a{\alpha} 
\def\t{\tau}
\def\b{\beta} 

\def\xap{x_\a^+}
\def\xam{x_\a^-}
\def\1{{\bf 1}}
\def\la{\langle}
\def\ra{\rangle}
\def\rtar{\rightarrow}
\def\mt{\mapsto} 
\def\pmatrix{\left(\begin{array}{cc}}
\def\endpmatrix{\end{array}\right)}

\def\half{\frac{1}{2}}
\def\fourth{\frac{1}{4}}
\def\fffourth{\frac{3}{4}}
\def\third{\frac{1}{3}}
\def\tthird{\frac{2}{3}}
\def\eighth{\frac{1}{8}}
\def\eeeighth{\frac{3}{8}}

\def\sixteenth{\frac{1}{16}}
\def\sssixteenth{\frac{3}{16}}  

\def\thirtysecond{\frac{1}{32}}
\def\ttthirtysecond{\frac{3}{32}}

\def\eop{{$\square$}}

\def\dual#1{#1^*} 

\def\dg#1{{\cal D}(#1)}

\def\kron#1#2{\delta_{#1#2}}

\def\avl{Aut(V_L)} 
\def\avlp{Aut(V_L^+)} 
\def\vlp{V_L^+}
\def\lvoap{LVOA${}^+$} 

\def\envalg#1{{\cal U}( #1 )}

\def\ve{\vfill \eject}

\centerline { 
21 September, 2004
}

\bigskip

{\bf Abstract.} Let $L$ be a positive definite even lattice 
and $V_L^+$ be the fixed points of the lattice VOA $V_L$ associated to 
$L$ under an  automorphism of $V_L$ lifting the $-1$ isometry of
$L$.   For any positive rank, 
the full automorphism group of $V_L^+$ is determined if
$L$ does not have vectors of norms 2 or 4.  For any $L$ of rank 2, a set of 
generators and the full automorphism group of $V_L^+$ are determined.

\ve
\tableofcontents 
\ve

\section{Introduction}

This article continues a program to study automorphism groups 
of vertex operator algebras.  See references in the survey \cite{G2} and
the more recent articles \cite{G1}, \cite{DG1}, \cite{DG2}, \cite{DGR}
and \cite{DN1}.

Here we investigate the fixed point 
subVOA of a lattice type VOA with respect to a group of order 
2 lifting the $-1$ map on a positive definite lattice.  We can obtain 
a definitive answer for the automorphism group
of this subVOA in two extreme cases.  The first is where the
lattice has no vectors of norms 2 or 4, and the second is 
where the lattice has rank 2.  

We use the standard notation $V_L$ for a lattice VOA, based on the
positive definite even integral lattice, $L$.  For a subgroup $G$ 
of $Aut(L)$, $V_L^G$ denotes the subVOA of points fixed by $G$.
When $G$ is a group of order 2 lifting $-1_L$, it is customary
to write $V_L^+$ for the fixed points (though, strictly speaking,
$G$ is defined only up to conjugacy; 
see the discussion in \cite{DGH} or \cite{GH}).

The rank 2 case is a natural extension of work on 
the rank 1 case, where $Aut(V_L^G)$ was determined for all 
rank 1 lattices $L$ and all choices of finite group $G\le Aut(V_L)$.  
 The styles of proofs are different.  In the rank 1 case, 
there was heavy 
analysis of the representation theory of the principal
Virasoro subVOA on the ambient VOA.  
In the rank 2 case, 
there is a lot of work on idempotents, solving 
nonlinear equations as 
well as work with several subVOAs associated to 
Virasoro elements.   
For rank 2, the case of nontrivial degree 1 part is harder to settle 
than in rank 1.

Our strategy follows this model.  Let $V$ be one of our $V_L^+$.
We get information about $G:=Aut(V)$ by its action on the finite 
dimensional algebra $A:=(V_2, 1^{st})$.  We take a subset $S$ of 
$A$ which is $G$-invariant and understand $S$ well enough to limit the 
possibilities for $G$ (usually, there are no automorphisms besides 
the ones naturally inherited from $V_L$).   A natural choice for $S$ is 
the set of idempotents or conformal vectors.  
Usually, $S$ spans $A$, or 
at least generates $A$.
In the main case of rank 2 lattice, we prove that 
$Aut(V)$ fixes a subalgebra of $A$ which is the natural $M(1)^+_2$.
The structure of $V$ is controlled by $M(1)^+$, which 
is generated by $M(1)^+_2$ and its eigenspaces, so we eventually
determine 
$G$.   

For several results, we give more than one proof.  

For the case of a lattice $L$ without roots, 
the automorphism group of $V_L^+$ was studied in the recent article 
\cite{S}.

We thank Harm Derksen for help with computer algebra.

\section{Background Definitions and Notations}

\begin{nota}\labttr{latticenota} 
Let $L$ be an even integral lattice.  For an integer $m$, define
$L_m:=\{ x \in L | (x,x)=2m\}$.  Let $H:=\CC \otimes L$, the ambient complex vector space.  
For a subset $S$ of $L$, define  $rank(S)$
 to be the rank of the sublattice
spanned by $S$.
\end{nota}

\begin{de}\labttr{autl} For a lattice, $L$, the group of automorphisms of 
the free abelian group $L$ which preserves the bilinear form is 
called {\it the group of automorphisms, the isometry group,  the group of units}
or {\it the orthogonal group of $L$}.  
This group is denoted $Aut(L)$ or  $O(L)$.  
We will use the notation $O(L)$ in this article, as well as the
associated $SO(L)$ for the elements of determinant 1, $PO(L)$ for 
$O(L)/\{\pm 1 \}$ and $PSO(L)$ for $SO(L)/SO(L) \cap \{\pm 1 \}$.  
\end{de}

\begin{de} \labttr{autlhat} For an even integral lattice,
$L$,  we let $\hat L$ be the 2-fold cover of $L$ described
in \cite{FLM},  \cite{DGH}, 
\cite{GH}.  We may write bars for the map $\hat L \rightarrow L$.  
The  {\it the group of automorphisms, the isometry group,  the group of units}
or {\it orthogonal group} is the set of group automorphisms 
of $\hat L$ which preserve the bilinear form on the quotient of $\hat L$ by
the normal subgroup of order 2. It is denoted $Aut(\hat L)$ or $O(\hat L)$ 
and has shape $2^{rank(L)}.O(L)$.  We use bars to denote the natural map $O(\hat L)\rightarrow O(L)$.  
\end{de}  

We list some notations for work with lattice type VOAs.  

\centerline{\bf List of Notations} 

\bigbreak
\halign{#\hfil&\quad#\hfil\cr

$\dg L$  &  the discriminant group
of the integral lattice $L$ is $\dg L :=\dual L / L$. \cr 

$e^\a$ & standard basis element for $\CC [L]$ \cr 

FVOA	& framed vertex operator algebra  \cite{DGH}
 \cr 
LVOA  & lattice vertex operator algebra  \cite{FLM}\cr

LVOA type	  &  the fixed points of a 
		lattice vertex operator algebra under a \cr & finite group 
		of automorphisms \cite{DG1, DGR}\cr
LVOA+ & $V_L^+$ for an even lattice $L$\cr

LVOAG$(L)$	&  the subgroup of $Aut(V_L)$, for an even integral 
		lattice $L$, \cr & as described in \cite{DN1};
		it is denoted $\NN (\hat L)$ and is
		 an extension\cr &  of the form $T.Aut(L)$ (possibly 
		nonsplit), where $T$ is a natural copy 
\cr & of the 
		torus $\CC \otimes L/\dual L$ 
\cr &
obtained by 
		exponentiating the maps $2\pi x_0$, for $x \in V_1$;
\cr &  the
		quotient of this group by the normal subgroup \cr & $T$ is
		naturally isomorphic to $Aut(L)$.  Also, $\NN (\hat L)$ is 
		the product of 
\cr & subgroups $TS$, where 
		$S\cong O(\hat L)$
		and 
\cr & $S \cap T = \{ x \in T| x^2=1\}\cong \ZZ_2^{rank(L)}$.  
		We may take $S$ to be 
		 the 
\cr & centralizer in LVOAG$(L)$ of a lift of $-1$; it
                has the form
\cr &  $2^{rank(L)}.Aut(L)$ and in fact any such 
		$S$ has this form.  
\cr & Denote the groups $S, T$ by $\OO (\hat L)$ and $\TT(\hat L)$,
           respectively.   
\cr 

LVOA group for $L$ &	this means LVOAG($L$).   \cr

LVOAG &		this means LVOAG$(L)$, for some $L$ \cr

LVOAG${}^+(L)$ &  this is the centralizer in LVOAG$(L)$ of a lift of $-1$ 
		modulo the
\cr &  group of order 2 generated by the lift; 
\cr & it
		has the form $2^{rank(L)}.[Aut(L) / \la -1 \ra ]$; it is the 
		{\it inherited group}  \cr

LVOAG${}^+$  &	this means LVOAG${}^+(L)$, for some $L$.   \cr

LVOA${}^+$-group &	same as  LVOAG${}^+$  \cr

$M(1), M(1)^+$ &  See Section 3.  \cr

$\NN (\hat L)$ & See LVOAG$(L)$ \cr

$o$ &		linear map from $V$ to $End(V)$ \cr 

$\OO (\hat L)$ & See LVOAG$(L)$ \cr

$\TT  (\hat L)$ & See LVOAG$(L)$ \cr

$v_{\alpha}$ & $e^{\a}+e^{-\a}$ \cr

$\XX$ or $\XX (L)$: & given an even integral lattice, $L$,
this is a group of shape $2^{1+rank(L)}$ 
\cr & for which
commutation corresponds to inner products modulo 2; 
\cr & see an
appendix of
\cite{GH}.  \cr

$\XX \OO$ or $\XX \OO (\hat L)$ & an extension of $\XX$ upwards by $O(L)$. \cr

$\XX \PP \OO$ or  $\XX \PP \OO (\hat L)$&  a quotient of $\XX \OO$ by a central
involution which corresponds to $-1_L$
\cr & under the natural 
epimorphism to $O(L)$. \cr
}

\begin{rem}\labttr{splitting} 
If $(L,L)\subset 2\ZZ,$ $\hat L  \cong L\times\<\pm 1\>.$ 
Thus $O(\hat L)$ contains a copy of  $O(L)$ which
complements the normal subgroup of order $2^{rank(L)}$
consisting of automorphisms which are trivial on the
quotient group $L$ of $\hat L$. 
This splitting passes to the groups $PO(\hat L)$ and $\XX
\PP \OO (L)$.  
\end{rem}

\section{Automorphism group of $V_L^+$ with $L_1=L_2=\emptyset$  }

In this section,  we determine the automorphism group of $V_L^+$ 
with $L_1=L_2=\emptyset$ and assume only
that $rank (L)>1.$ The automorphism
group of $V_L^+$ in the case $rank (L)=1$ is determined in [DG1]
without any restriction  on $L.$ The assumption that
$L_1=L_2=\emptyset$ ensures that any automorphism of $V_L^+$
preserves the subspace $M(1)^+_2$,   
which can be identified
with the Jordan algebra $S^2H$.  

Since $M(1)^+$ is generated by
$M(1)_2^+$ if
$\dim H>1$ and $V_L^+$ is a direct sum 
of eigenspaces for $M(1)_2^+$ (cf. [AD]), the structure of $Aut (V_L^+)$ 
can be determined easily. 
We shall use a classic result.

\begin{prop}\labtt{autsymmat} The automorphism group of the
Jordan algebra of  symmetric $n\times n$ matrices is
$PO(n,\CC )$, acting by conjugation.
\end{prop}
\pf  \cite{J}.  
\eop

\subsection{$Aut (M(1)^+)$}

We first recall the construction of $M(1)^+.$ Let $H$ be a 
$n$-dimensional complex vector 
space with a nondegenerate symmetric bilinear form
$(\cdot , \cdot)$ and $\hat H=H\otimes 
\CC[t,t^{-1}]\oplus \CC\,c$ the corresponding affine Lie algebra.
Consider the
induced $\hat H$-module  
$$M(1)= 
\envalg{\hat H}\otimes_{\envalg{H\otimes{{\CC}}[t]
\oplus{{\CC}}c}}{{\CC}}\simeq S(H\otimes t^{-1}\CC[t^{-1}])\ \ \ (\mbox{linearly})
$$
where $H\otimes {{\CC}}[t]$ acts trivially on $\CC,$
and $c$ acts as 1. For $\alpha\in H$ and $n\in\ZZ$ we set 
$\alpha(n):=\alpha\otimes t^n.$ Let $\tau$ be the
automorphism of $M(1)$ such that
$$\tau ( \a_1(-n_1)\cdots \a_k(-n_k))=(-1)^k\a_1(-n_1)\cdots \a_k(-n_k)$$
for $\a_i\in H$ and $n_1\geq \cdots \geq n_k\geq 1.$ Then $M(1)^+$ 
is the fixed point subspace of $\tau.$

\begin{prop}\labtt{pm1} The automorphism  group of $M(1)^+$ is $PO(n,\CC).$
\end{prop}
\pf We first deal with the case that $\dim H>1.$ Then 
$M(1)^+$ is generated by $M(1)^+_2$ (cf. [DN2]),  which is a
Jordan algebra under $u\cdot v=u_{1}v$ for $u,v\in
M(1)^+_2.$ So any automorphism of $M(1)^+$ restricts to an
automorphism of the Jordan algebra
$M(1)^+_2.$ On the other hand, the automorphism group of $M(1)$
is $O(n,\CC)$ [DM2], which preserves $M(1)^+.$ Clearly, the kernel of
the action of $O(n,\CC)$ on $M(1)^+$ is $\{ \pm 1\}.$ As a result $PO(n,\CC)$
is a subgroup of the automorphism group of $M(1)^+.$  By Proposition
\ref{autsymmat}, 
any automorphism of $M(1)^+_2$ extends to an automorphism
of $M(1)^+.$ 

We now assume that $\dim H=1.$ Then $M(1)^+$ is not generated by $M(1)^+_2.$
By Lemma 2.6 and Theorem 2.7 of [DG1] for any nonnegative even integer
$n$ there is a unique lowest weight vector $u^n$ (up to scalar multiple)
of weight $n^2$ and $M(1)^+$ is generated by the Virasoro vector and
$u^n.$   Using the fusion rule given in Lemma 2.6 of
[DG1] we immediately  see that the automorphism group of $M(1)^+$ in this case is
trivial. Clearly, $PO(1,\CC)=1.$
This finishes the proof.
\eop

\subsection{Aut$(V_L^+)$ }

First we review from [B] and [FLM] the construction of lattice vertex operator algebra $V_L$
for any positive definite even lattice $L.$ Let $H=\CC\otimes_{\ZZ}L$.  
Recall that $\hat{L}$ is the canonical central 
extension of $L$ by the cyclic
group $\< \pm 1\>$ such that the commutator map is given 
by $c(\alpha,\beta)=(-1)^{(\a,\b)}.$ We fix a bimultiplicative
2-cocycle $\epsilon: L\times L\to \<\pm 1\>$ such that
$\epsilon(\a,\b)\epsilon(\b,\a)=c(\a,\b)$ for $\a,\b\in L.$ 
Form the induced
$\hat{L}$-module
\begin{eqnarray*}
{\CC}\{L\}=\CC[\hat{L}]\otimes _{\CC [ \< \pm 1\> ]}\CC\simeq
\CC[L]\;\;\mbox{(linearly)},\end{eqnarray*}
where $\CC[\cdot]$ denotes the group algebra and $-1$ 
acts on $\CC$ as multiplication by $-1$. 
For $a\in \hat{L}$, write $\iota (a)$ for
$a\otimes 1$ in $\CC\{L\}$. 
Then the action of $\hat{L}$ on ${\CC}
\{L\}$ is given by: 
$a\cdot \iota (b)=\iota (ab)$ 
for $a,b\in \hat{L}$.  
If $(L,L)\subset 2\ZZ$ then 
$\CC\{L\}$ and $\CC[L]$ are isomorphic algebras. 
The lattice vertex operator algebra $V_L$ is defined to be 
$M(1)\otimes
\CC\{L\}$, as a vector space. 
 
Then $O(\hat L)$ is a naturally defined  
subgroup of $\Aut(\hat L)$ and
$Hom(L,\ZZ/2\ZZ)$ may be identified with  a subgroup of $O(\hat L)$ (see [FLM] , [DN1], [GH]) and
there is 
 an exact sequence
$$1\rightarrow Hom (L,\ZZ/2\ZZ) \rightarrow
O(\hat{L})\overset{-}{\rightarrow}
O(L)\rightarrow 1.$$
It is proved in [DN1] that $Aut (V_L)$ has shape $N\cdot O(\hat L)$   
where
$N$ is the normal subgroup of $Aut (V_L)$ generated by $e^{u_0}$ for
$u\in (V_L)_1.$  Note that $Hom (L,\ZZ/2\ZZ) $ can furthermore be identified with the
intersection of $N$ and $O(\hat L)$.  See the List of Notations.  

Let $e:L \rightarrow \hat{L}$ be a section associated to the 2-cocycle 
$\epsilon$, written $\a \mapsto e_\a$. Let $\theta$ be the automorphism $\theta$
of $\hat L$ of order $2$ such that $\theta e_{\a}=e_{-\a}$ for
$\alpha \in L.$ Then $\theta$ extends to an automorphism
of $V_L,$ still denoted by $\theta$, such that 
$\theta|_{M(1)}$ is identified with $\tau$  and
$\theta \iota(a)=\iota(\theta a)$ for all $a \in \hat L$.   
Set $e^{\alpha}=\iota(e_{\a}).$ Then $\theta e^{\alpha}=e^{-\a}.$ 

Let $V_L^+$ be the fixed points of $\theta.$ In order to determine
the automorphism group of $V_L^+,$ it is important to understand 
which automorphism of $V_L$ restricts to an automorphism of $V_L^+.$
Clearly, the centralizer of $\theta$ in $Aut (V_L)$ acts on $V_L^+.$
So we get an action of  $O(\hat L)/\<\pm 1\>$ on $V_L^+.$ 
Let $h\in H.$ Then $e^{2\pi ih(0)}$ preserves $V_L^+$ if and only if
$(h,\alpha) \equiv (h,-\alpha)$ modulo $\ZZ$ for any $\alpha\in L.$ That is,
$h\in \frac{1}{2}\dual L$ where $\dual L$ is the dual lattice of
$L.$    

\begin{lem}\labtt{lema} 
The subgroup of $\avlp$ which preserves 
$M(1)^+_2$ is just the LVOA${}^+$-group.  
\end{lem}

\pf  Let $n:=dim (H)$.   Let $\sigma\in Aut(V_L^+)$  such that $\sigma M(1)^+_2\subset M(1)^+_2.$
 Then $\sigma|_{M(1)_2^+}\in PO(n,\CC)$ as in  \ref{pm1}. Note that $M(1)^+$ is generated by $M(1)_2^+$ as  $rank(L)>1$ (see the proof
of Proposition \ref{pm1}). So, $\sigma$ preserves $M(1)^+.$
 
For any $\alpha\in L$, let $V_L^+(\alpha)$ be the $M(1)^+$-submodule generated by
$v_{\alpha}:=e^{\a}+e^{-\a}.$ Then $V_L^+(\alpha)$ is an irreducible
$M(1)^+$-module, $V_L^+(\alpha)$ and $V_L^+(\b)$ are isomorphic $M(1)^+$-modules
if and only if $\alpha=\pm \beta$ (cf. \cite{AD}).
Moreover, if $\alpha\ne 0$ then $V_L^+(\alpha)$ is isomorphic to 
$M(1)\otimes e^{\alpha}$ (cf. \cite{AD}). 

Note that $V_L^+=\sum_{\alpha\in L}V_L^+(\alpha)$.  
Let $S$ be a subset
of $L$ such that $|S\cap \{\pm\alpha\}|=1$ for any $\alpha\in L.$ Then
for any two different $\alpha,\beta\in S,$ $V_L^+(\a)$ and $V_L^+(\b)$
are nonisomorphic $M(1)^+$-modules and   
$$V_L^+=\oplus_{\alpha\in S}V_L^+(\alpha)$$
is a direct sum of nonisomorphic irreducible $M(1)^+$-modules. 

Let $\alpha\in L.$ Since $\sigma$ preserves $M(1)^+,$ it sends $V_L^+(\a)$ to 
$V_L^+(\beta)$ for  some $\beta\in L.$ The vector  $v_{\alpha}$ is the 
unique lowest weight vector (up to a scalar) of $V_L^+(\alpha).$  
This implies that  $\sigma(v_\alpha)=\l v_{\beta}$ for some nonzero 
scalar $\l \in \CC$  (depending on $\a$ and $\b$).  

For a vertex operator algebra $V$ and a homogeneous
$v\in V$ we set $o(v)=v_{\wt v-1}$ and extend to all of $V$ linearly.
Note that  $v_{\alpha}$ is an eigenvector for 
$o(v)$ for $v\in M(1)^+_2.$  In fact,
$o(h_1(-1)h_2(-1))v_{\alpha} =(h_1,\alpha)(h_2,\alpha)v_{\alpha}$ 
for $h_i\in H.$ 
Recall the proof of Proposition \ref{pm1}. 
We can regard the restriction of 
$\sigma$ 
to $(V_L^+)_2\cong  M(1)^+_2$ 
as an element of $O(n,\CC)$, well-defined modulo $\pm 1$.   
Then $\sigma(h_1(-1)h_2(-1))=
(\sigma h_1)(-1)(\sigma h_2)(-1).$  Note that $\s^{-1}$ is the adjoint of $\s$.  Then, 
\begin{eqnarray*}
& &(h_1,\alpha)(h_2,\alpha)\lambda  v_{\beta}= \s ((h_1,\alpha)(h_2,\alpha)  v_{\a}) = 
\sigma ( o(h_1(-1)h_2(-1))v_{\alpha})\\
& & =o(\sigma(h_1(-1)h_2(-1)))\lambda v_\b =(\sigma h_1,\beta)(\sigma h_2,\b)\lambda v_\b .\\
\end{eqnarray*}
Since the $h_i$ are arbitrary,  $\sigma\a=\pm \beta.$ 
Thus $\sigma$ maps $L$ onto $L$ so induces an isometry of $L$ which is well defined  modulo $\<\pm 1\>.$ 

Multiplying $\sigma$ by an element from LVOAG${}^+(L)$ (which comes from $\NN (\hat L )$), 
we can assume that
$\sigma|_{M(1)^+}=id_{M(1)^+}$ Then $\sigma v_\a=\lambda_{\a} v_\a$ for
some nonzero $\lambda_\a\in \CC.$ Since $V_L^+(\alpha)$ is an
irreducible $M(1)^+$-module we see that $\sigma$ acts as the scalar $\lambda_\alpha$
on $V_L^+(\alpha).$ Clearly, $\l_\a=\l_{-\a}.$ 
Note that
\begin{eqnarray*}
& &Y(v_{\alpha},z)v_{\beta}=E^-(-\a,z)\epsilon(\alpha,\beta)e^{\a+\b}z^{(\a,\b)}+
E^-(-\a,z)\epsilon(\alpha,-\beta)e^{\a-\b}z^{-(\a,\b)}\\
& &\ \ +E^-(\a,z)\epsilon(-\alpha,\beta)e^{-\a+\b}z^{-(\a,\b)}+
E^-(\a,z)\epsilon(\alpha,\beta)e^{-\a-\b}z^{(\a,\b)}
\end{eqnarray*}
where 
$$E^-(\alpha,z)=\exp(\sum_{n<0}\frac{\alpha(n)z^{-n}}{n}).$$
Thus, if $n$ is sufficiently negative,    
$(v_\a)_n(v_\b)=u+v$ for some nonzero $u\in V_L^+(\a+\b)$
and $v\in V_L^+(-\a+\b).$   
This gives 
$\l_\a\l_\b=\l_{\a+\b}=\l_{\a-\b}$ by applying
$\sigma$ to $(v_\a)_n(v_\b)=u+v.$ 
So $\alpha\mapsto \l_\a$ defines a character of
abelian group $L/2L$ of order $2^n.$ Clearly, any character
$\l: L/2L\to \<\pm 1\>$ defines an automorphism $\sigma$ 
which acts on $V_L^+(\a)$ as $\l_\alpha.$
As a result,
the subgroup of $Aut (V_L^+)$ which acts trivially on
$M(1)^+$ is isomorphic the dual group of $L/2L$ and 
is exactly the subgroup  of $O(\hat L)/\<\pm 1\>$ which we identified as 
$\Hom(L,\ZZ/2\ZZ)$.     
As a result, the subgroup of $Aut (V_L^+)$ which preserves
$M(1)^+_2$ is exactly the group $O(\hat L)/\<\pm 1\>,$ as desired.
\eop

\begin{prop}\labtt{hrank} Let $L$ be a positive definite even lattice such
that $L_1=L_2=\emptyset$.  Then $Aut (V_L^+)$ is the inherited group, i.e., the 
LVOA${}^+$-group.  
\end{prop}

\pf In this case we have $(V_L^+)_2=M(1)^+_2.$ Thus, any automorphism
of $V_L^+$ preserves $M(1)^+_2.$ By Lemma \ref{lema}, $Aut (V_L^+)$ is
the LVOA${}^+$-group.
\eop

\section{Rank 2 lattices}

All lattices in this article 
are positive definite.  
Throughout  this article, $L$ denotes an even integral lattice.  
We recall a general result.

\begin{lem}\labtt{detlattice}  Let $L$ be a lattice and $M$ a
sublattice.

(i)  If $|L:M|$ is finite, $det(M)=det(L)|L:M|^2$.

(ii)  If $M$ is a direct summand of $L$, $L/[M+ann_L(M)]$ embeds in 
$\dg M$.
\end{lem}  

\pf  These are standard results.  For example, see \cite{G3}.
\eop

We need to sort out rank 2 lattices by whether they contain roots or elements 
of order 4, due to their contributions to low degree terms of the lattice
VOA.   
We shall use the notations \ref{latticenota}.

\begin{lem}\labttr{rank2}
Suppose that $rank(L_1)=2$.  Then $L_1$ spans $L$ and $L$ is one of
$L_{A_1^2}$ or $L_{A_2}$.
\end{lem}
\pf The span of $L_1$ is isometric to $L_{A_1^2}$ or $L_{A_2}$.  Each of these is
a maximal even integral lattice under containment.
\eop 

\begin{lem}\labttr{rank1} Suppose that $rank(L_1)=1$.  Let $r \in
L_1$ and let 
$s$ generate $ann_L(r)$.  
Then $(s,s)\ge 4$ and if $L>  span\{r, s\}$, then $14 \le (s,s) \in 6+8\ZZ$.    
\end{lem}
\pf 
Note that $\ZZ r$ is a direct summand of $L$.  
We have $(s,s) \ge 4$.  In case $L > N:= span\{r, s\}$, $L/N$ has order 2, by
\ref{detlattice}.  
If $x$ represents the nontrivial coset, $(x,x)\ge 4$ then
$(2x,2x)\ge 16$.  Also, $(2x,2x)\in 8\ZZ$.  
Since $(x,r)$ is odd, if we write $2x=pr+qs$, for integers $p, q\in \ZZ$, then $p$
is odd and so $q^2(s,s) \in 6+8\ZZ$. It follows that $q$ is odd and $(s,s)\in
6+8\ZZ$.  
\eop

\begin{lem}\labtt{rank02} Suppose that $L_1=\emptyset$ and
$rank(L_2)=2$.   If $r, s$ are linearly independent norm 4
elements, then  they span $L$ and have Gram matrix $G=\pmatrix  4&b\cr
b&4\endpmatrix$, for some $b \in \{0, \pm 1, \pm 2\}$.  
\end{lem}
\pf 
If $L\ne N:=span \{r,s\}$, then $det(N)=16-b^2$ is divisible by a perfect
square, whence 
$b=0$ or $b=\pm 2$ and the index is 2.  
Actually, $b=0$ does not occur here since $\half r, \half s \not \in L$ implies that $\half (r+s) \in L_1$, a
contradiction.  So, $b=\pm 2$.  
 Clearly, $span \{r,s\}\cong \sqrt 2 L_{A_2}$.  However, any 
integral lattice containing the latter with index 2 is odd, a 
contradiction.  Therefore, $L=N$ and the Gram matrix is as above.  Positive definiteness implies that $|b|< 4$ and rootlessness implies that $b\ne \pm 3$.  
\eop

\begin{lem}\labtt{rank01}  Suppose that $L_1=\emptyset$ and
$rank(L_2)=1$. Let $r\in L_2$ and let $s$ generate
$ann_L(x)$.  Then 
$(s,s) \ge 6$ and $L/ span \{r,s\}$ is a subgroup of $\ZZ_4$.

If the order of $L/ span \{r,s\}$ is $2, 8 \le (s,s)\in 4+8\ZZ$.

If the order of $L/ span \{r,s\}$ is $4, 28 \le (s,s)\in 28+32\ZZ$.  
\end{lem}
\pf
Let $x$ be in a nontrivial coset of $N:=span \{r,s\}$ in $L$.  

If $(x,r)\in
2+4\ZZ$, $2x=pr+qs$, where $p$ is odd.  We have $(x,x) \ge 6$, $p$ odd and $q\ne
0$.  Therefore, $(2x,2x) \in 8\ZZ$, $24\le 4p^2+(s,s)q^2$, whence $q$ is odd and
$(s,s)\in 4+8\ZZ$. 

If $(x,r) \in 1+2\ZZ$, $(4x,4x) \in 32\ZZ$.  We have $(x,x)\ge 6$, whence
$(4x,4x)\ge 96$.   If we write
$4x=pr+qs$, we have $4p^2+q^2(s,s) \in 32\ZZ$.  
Since $p$ is odd, 
$p^2\in 1+8\ZZ$
and $4p^2\in 4+32\ZZ$.  Since $(s,s)$ is even, $q$ is odd, $q^2\in 1+8\ZZ$ and
$(s,s) \in 4+8\ZZ$.  It follows that $\fourth q^2(s,s) \in 7+8\ZZ$ whence $\fourth
(s,s) \in 7+8\ZZ$.  
\eop

\section{About idempotents in small dimensional algebras}

We can derive a lot of information about the automorphism group of a vertex
operator algebra by restricting to low degree homogeneous pieces.  For the $V_L^+$
problem, the degree 2 piece and its product $x, y \mapsto x_1y$ give an algebra
which is useful to study.  Here,  for $rank(L)=2$, we concentrate on some
commutative algebras of dimension around 5.  Commutativity of $(V_L^+, 1^{st})$ is implied if $L_1=\emptyset$, which is so for $b \in \{ 0, \pm 1, \pm 2\}$ as in \ref{rank02}.  

It does not seem advantageous to give particular values to $b$ most of the time,
so we keep it as an unspecified constant in case these arguments might be a model for future work.  In the present work, we shall note limits on $b$, as needed.

\begin{nota}\labttr{5dimnota} Let $S$ be the Jordan algebra
of degree 2 symmetric matrices and suppose that $A$ is a
commutative 5 dimensional algebra of the form 
$A=S \oplus \CC v_r \oplus \CC v_s$.  Suppose that $v_r\times v_s=0$ and that  the notations of Appendix: Algebraic rules 
apply here, with the usual inner products and algebra product.  
Let $w=p+c_rv_r+c_sv_s$ be an idempotent.   
Suppose also that
$t$ is a norm  4 vector orthogonal to $r$.  
Let $a_1, a_2, a_3$ be scalars so that 
$p=a_1r^2+a_2rt+a_3 t^2$.    
\end{nota}

\begin{rem}\labttr{identityrs}
We note that the basis $r, s$ of $H$ has dual basis 
$r^*, s^*$, where $r^* = \frac{1}{16-b^2} (4r-bs)$ 
and
$s^*=\frac{1}{16-b^2} (4s-br)$.
The identity of $A$ is $\fourth \frac{1}{16-b^2} (rr^*+ss^*)
= \fourth \frac{1}{16-b^2} (4r^2 + 4s^2 -2b rs)$.  
\end{rem}

\begin{nota}\labttr{pq}
If $w$ is an element of $A$, write $w=p+q$ for $p \in S^2H$
and
$q\in \CC v_r \oplus \CC v_s$.  
Call the element $\bar w := p - q$ the {\it  conjugate
element}.  The components $p, q$ are called the {\it $P$-part} 
and the {\it $Q$-part}  of $w$.  
Extend this notation to subscripted elements: 
$w_i=p_i+q_i$, $\bar {w_i}=p_i-q_i$, for indices $i$.  
\end{nota}  

\begin{rem}\labttr{pq2}   In \ref{pq}, $q^2\in S^2H$ since $v_r\times v_s=0$.  Also $w=p+q$ is an idempotent if and only if $p=p^2+q^2$ and
$q=2p\times q$.   Therefore, $w=p+q$ is an idempotent if and only if the conjugate $p-q$
is an idempotent.  
\end{rem}

\begin{lem}\labttr{sumofidemps}  Suppose that
$w_1$ and $w_2$ are idempotents 
and their sum is an idempotent.  Then $w_1\times w_2=0$ and
$(w_1,w_2)=0$.  
\end{lem} 
\pf 
We have $(w_1+w_2)^2=w_1^2+2w_1\times w_2 + w_2^2$, whence
$w_1\times w_2=0$.  Also,
$(w_1,w_2)=(w_1^2,w_2)=(w_1,w_1\times w_2)=0$.  
\eop

\begin{de}\labttr{types012} Throughout this article, an
idempotent is not zero or the identity, unless the context
clearly allows the possibility.   We call an idempotent
$w$ of {\it type 0, 1, 2}, respectively, if it has $Q$-part which is 
0, is a multiple of $v_r$ or $v_s$, or is not a multiple of either  $v_r$
or
$v_s$.
\end{de}

\begin{lem}\labttr{ad=4}  Then 
(i) $r^2\times
s^2=4brs$, $r^2\times r^2 = 16r^2$, $s^2\times s^2=16s^2$,
$rs\times rs=4r^2+4s^2+2brs$;  
$x^2 \times v_r= (x,r)^2v_r=\half (x^2,r^2)v_r$; 
$r^2\times rs= 8rs+2b^2 r^2$, $s^2\times rs=8rs+2b^2 s^2$;  
also
$v_r\times v_r=r^2, v_r\times v_s=0, v_s\times v_s = s^2$.

(ii) $(r^2,r^2)=32=(s^2,s^2)$, $(r^2,s^2)=2b^2$, $(rs,
rs)=16+b^2$, 
$(rs, r^2)=8b=(rs, s^2)$ and   $(v_r,v_r)=2=(v_s,v_s)$ and
$(v_r,v_s)=0$.  
\end{lem} 
\pf See the Appendix (and take $a=d=4$, $b\ne 0, \pm 2$).
\eop

\subsection{Idempotents of type 0}

\begin{lem}\labttr{type0}
These are just idempotents in the Jordan algebra of
symmetric matrices.  They are ordinary idempotent matrices
which are symmetric.  Up to conjugacy by orthogonal 
transformation, they are diagonal matrices with diagonal
entries only 1 and 0.  
\end{lem}

\subsection{Idempotents of type 1}

\begin{nota}\rm The next few results apply to the case of an idempotent of type 1, i.e., 
the form $w=p+c_rv_r$, where $c_r\ne 0$.  In such a case,
$w=w^2=p^2+c_r^2r^2+c_r(p,r^2)v_r$ (see the Appendix : Algebraic Rules)  From $c_r\ne 0$, we get
$(p,r^2)=1$.  We continue to use the notation of \ref{5dimnota}.  
\end{nota}

\begin{lem} \labtt{a123} 
Suppose that $c_r\ne 0$ and $c_s=0$.   
We have 
$a_1=16 a_1^2 + 4a_2^2+c_r^2$;  $a_2=16a_2(a_1+a_3)$; 
$a_3=16 a_3^2 + 4a_2^2$ and $(p,r^2)=1$.
\end{lem}
\pf
Compute $p+c_rv_r=w=w^2= 
p^2+c_r^2r^2+c_r(p,r^2)v_r$ (see \ref{ad=4}) and expand in the basis 
$r^2, rt, t^2, v_r$.  
\eop 

\begin{coro} \labtt{a1}$a_1={1\over 32}$.
\end{coro}
\pf 
We have $1= (p,r^2)=a_1(r^2,r^2)=32a_1$, whence $a_1={1\over 32}$.
\eop

\begin{lem} \labtt{a123b} 
Suppose that $c_r\ne 0$ and $c_s=0$.   Then 

(A1) 
$a_1=16a_1^2+4a_2^2+c_r^2$; 

(A2) $a_2=16a_2(a_1+a_3)$; and 

(A3) 
$a_3=16a_3^2+4a_2^2$.  
\end{lem} 
\pf
Compute
$p^2=(16a_1^2+4a_2^2)r^2+16(a_1a_2+a_3a_2)rt+(16a_3^2+4a_2^2)t^2$
and use $w=w^2=p^2+c_r^2 r^2 +c_rv_r$.
\eop

\begin{lem}  \labtt{a2=0}
Suppose that $c_r\ne 0$ and $c_s=0$.   If $a_2=0$, then 
$a_3\in \{ 0, \sixteenth \}$ and $c_r=\pm \eighth$.  
\end{lem}
\pf
We deduce from (A3) that $a_3=16a_3^2$, then $c_r=\pm
\eighth$.  
\eop

\begin{lem} \labtt{a2not0} 
Suppose that $c_r\ne 0$ and $c_s=0$.   Then $a_2=0$.
 \end{lem}  
\pf 
  If $a_2\ne 0$, then from (A2), $1=16(a_1+a_3)$ and  we get $a_3=\thirtysecond$. 
Next, use (A3) to get $\frac{1}{64}=4a_2^2$.  Finally use
(A1) to get $c_r=0$.  
\eop

\begin{thm} \labtt{maintype1} 
Assume that $c_r\ne 0$ and $c_s=0$.  
Then 

(i) $a_1=\thirtysecond, a_2=0, c_r=\pm \eighth$; and 

(ii) either $(w,w)=\sixteenth$ and $a_3=0$; or $(w,w)={3\over 16}$ and $a_3=\sixteenth$.  

All of the above cases occur.   If an idempotent occurs, so does its 
complementary idempotent.  
\end{thm} 
\pf  This is a summary of preceding  results.  
\eop

\begin{lem}\labtt{eigenspacestype1}
If $w$ is an idempotent of type 1, then 

(i) if $(w,w)=\sixteenth$, 
the eigenvalues for $ad(w)$ are $1,  0, 0,  \fourth, {\frac {b^2}{32}}$.  
Eigenvectors for these respective eigenspaces are 
$w, 1-w, t^2, rt,  v_s$;  

(ii) if $(w,w)=\sssixteenth$, 
the eigenvalues for $ad(w)$ are $0, 1,  1,  \fffourth, {1-\frac {b^2}{32}}$. 
Eigenvectors for these respective eigenspaces are 
$w, 1-w, t^2, rt,  v_s$.  

If  $ {\frac {b^2}{32}}\ne 0, 1, \fourth, \fffourth$, the multiplicities of 0 and 1
are 2 and 1 in case (i) and 1 and 2 in case (ii). 
\end{lem}
\pf
Straightforward calculation.  Note that   $ {\frac {b^2}{32}}\ne 0, 1, \fourth, \fffourth$ follows if $b$ is rational 
\eop

\begin{coro} \labtt{summandtype1}  
If $w$ is a type 1 idempotent and  is  the sum of two nonzero idempotents, 
$w_1, w_2$, then
$w$ has the form 
$\thirtysecond r^2 +\sixteenth t^2 \pm \eighth v_r$ 
and $w_1, w_2$ are, up to order,  
$\thirtysecond r^2 \pm \eighth v_r$ 
and $\sixteenth t^2$.  
\end{coro}
\pf
If $w$ is such a sum, 
each $w_i$ is in the 1-eigenspace of $ad(w)$, which must be more than 
1-dimensional.  This means that $w$ has norm $\sssixteenth$ and one of the $w_i$,
say for $i=1$,  has type 1 and $Q$-part $\pm \eighth v_r$.  Therefore, $w_2$ has
type 0, whence norm $\eighth$.  This means that $w_1$ has norm $\sixteenth$ and so
we know that $w_1$ has shape $\thirtysecond r^2 \pm \eighth v_r$ and $w=\sixteenth
t^2$.  
\eop


\subsection{Idempotents of type 2} 

\begin{hyp}\labttr{b=pm1} 
We assume in this subsection that the parameter $b \ne 0, \pm 2, \pm 3$ (which means $b=\pm 1$).  Then 
the algebra $(V_2, 1^{st})$ is commutative since $V_1=0$.    
\end{hyp} 

\begin{nota}\labttr{candd} 
$p=c(r^2+s^2)+d rs$, $v=c_rv_r+c_sv_s$.  
\end{nota}

\begin{lem}\labtt{crcsnonzero}  If $c_r$ and $c_s$ are
nonzero, then  there are at most 8 possibilities for $w$.   
In more detail, there are at most two values of $c$ 
(and, correspondingly, of $d$).  
We have $c_r^2=c_s^2$ and this common value depends on  $c$
(or on $d$).  
\end{lem}
\pf
Compute $p+c_rv_r+c_sv_s=w=w^2= 
p^2+c_r^2r^2+c_s^2s^2+c_r(p,r^2)v_r+c_s(p,s^2)v_s$.  
Since $c_r$ and $c_s$ are nonzero, 
$(p,r^2)= 1 =(p,s^2)$.  

Since $(r^2,r^2)=32=(s^2,s^2)$, $(rs, r^2)=8b=(rs, s^2)$ and 
$(r^2,s^2)=16+b^2$, we have $p=c(r^2+s^2)+d rs$.  
for some scalars, $c, d$.  The previous paragraph then
implies that $1 = (32+2b^2)c+8bd$.  
Since $b\ne 0$, 
$$d=\frac{1}{8b}(2c(32+2b^2)-1)\leqno (e1)$$ 
is a linear
expression in
$c$.  

Now, $p^2= (16c^2 + 4d^2)(r^2+s^2)+ (2bc^2+2bd^2 )rs$ and
so 
$w^2= (16c^2 + 4d^2+c_r^2) r^2+ (16c^2 + 4d^2+c_s^2)
s^2+ (2bc^2+2bd^2 )rs + 
c_rv_r+c_sv_s$.  It follows that $c_r^2=c_s^2$.

We compare coefficients of $r^2$ and get 
$$c=16c^2 + 4d^2+c_r^2. \leqno (e2)$$

We compare coefficients of $rs$ and get 
$$d=2bc^2+2bd^2.  \leqno (e3)$$  
Since $d$ is a linear expression in $c$, $c$ satisfies a
quadratic equation, depending on $b$ but not $c_r^2=c_s^2$. 
The degree of this equation really is  2 since 
$b\ne 0$ real implies that the  top coefficient  is nonzero.  

It follows that the ordered 
pair $d, c$ has at most two possible
values.  For each, there is a unique value for $c_r^2$,
hence at most two possible values for $c_r$ (and the same two
for $c_s$).  Therefore there are at most eight 
idempotents of type 2.  
\eop

\begin{lem}\labtt{cnot0} $c\ne 0$ and $d\ne 0$.  
\end{lem}
\pf  
Suppose that $c=0$.  
We then have $p=drs$, $d={\frac{-1}{8b}}$.  On the other hand, since $w=drs+v$ is
an idempotent, the coefficient for $w^2$ at $rs$ is 
$d=8bc^2+2bd^2=2bd^2=2bd^2$, which 
implies that $1=2bd$.  This  
is incompatible with  $d={\frac{-1}{8b}}$.  

If $d=0$, equation (e3) implies that $c=0$, which is false. 
\eop


\begin{lem}\labtt{coeffsofcomp}  
 If $w$ is a type 2
idempotent,
$w=c(r^2+s^2)+drs+c_rv_r+c_sv_s$ and $1-w$ is 
the complementary idempotent,
expanded similarly as 
$1-w=c'(r^2+s^2)+d'rs+c_r'v_r+c_s'v_s$, then 
$c_r'=-c_r\ne c_r$, $c_s'=-c_s$, 
$c\ne c'$ and $d\ne d'$.  
In particular, in the notation of \ref{crcsnonzero}, 
the function $c\mapsto c_r^2$ is
2-to-1 and so only one value of $c_r^2$ occurs 
for type 2 idempotents. 
\end{lem}
\pf
If it were true that $c=c'$, then $w=\half \II +v$ 
and $1-w=\half \II -v$.  
Since these are idempotents,  $v^2=\half \II$.  
However, this is
impossible as $b\ne \pm 2$ implies that $v^2$ 
is a multiple of $r^2+s^2$ and $\II$
is not a linear combination of $r^2, s^2$ for 
$b\ne 0$ (see \ref{identityrs}).  
\eop

\subsection{Sums of idempotents}

\begin{hyp}\labttr{stillb=pm1} 
We  continue to take $b=\pm 1$.  Results of the previous subsection apply.
\end{hyp} 

In the arguments in this section, we allow 
the symbol $b$ to be any odd integer, 
though the lattice is positive definite only for $b=\pm 1$.

\begin{lem}\labtt{sumoftwo}
Suppose that  $w_1, w_2$ are two idempotents of type 1.   If $w_1+w_2$ is an idempotent, then
$w_1+w_2$ does not have type 1 or type 2.   
\end{lem}
\pf 
We eliminate the sum having type 1 with Lemmas \ref{summandtype1}.  
To eliminate a sum having type 2, we note that for type 1 idempotents, 
we have $a_2=0$ by 
\ref{maintype1}, whereas $d\ne 0$ for type 2, by \ref{cnot0}.  
\eop

\begin{lem}\labtt{2+2notidemp} 
If $w_1, w_2$ are idempotents of type 2 and not
complementary, their sum is not an idempotent.
\end{lem}
\pf  Assume that the sum $w$ is an idempotent.  
From \ref{sumoftwo}, the sum has type 0, 
so has the form $\sixteenth u^2$, 
for some
vector 
$u \in H$ of norm 4.  
The eigenvalues of $ad(u)$ are $1,0,\half$ and
$\sixteenth (u,r)^2$,
$\sixteenth (u,s)^2$, with respective eigenvectors 
$u^2, 1-u^2, \half uu', v_r, v_s$,
where 
$u'$ spans the orthogonal of $u$ in $H$.  

Now, $w_1, w_2$ are  linearly independent (or else they
are equal, which is impossible).   This means that the eigenvalue 1 has
multiplicity at least 2.  So, at least one of
$(u,r)^2, (u,s)^2$ is 16.  Since $w_1, w_2$ lie in the 
1-eigenspace of $ad(u)$ and both $w_i$ have type 2, 
both these square norms must be 16, i.e., 
$m:=(u,r)=\pm 4$ and $n:=(u,s)=\pm 4$.     
Since $r, s$ form a basis and the form
is nonsingular, this forces $u=mr^*+ns^*$, 
where $r^*, s^*$ is the dual basis.
We have $4=(u,u)=16(r^*,r^*)+2mn(r^*,s^*)+16(s^*,s^*)$.  
The right side is 
$\frac{1}{16-b^2} [ 16(4r-bs,4r-bs)+2mn(4r-bs,4s-br)+16(4s-br,4s-br) ]$.  

Since $b$ is an odd integer, the above rational number in
reduced form clearly has numerator divisible  by 16, so does
not equal 4, a contradiction.
\eop

\begin{lem}\labtt{1+2notidemp}
The sum of a type 1 and type 2 idempotent is not an
idempotent.
\end{lem}
\pf 
Assume that $w:=w_1+w_2$ is an idempotent.  Obviously it does not have type 0.  
By
\ref{summandtype1}, it does not have type 1.  

We conclude that $w$ has type 2.  However, the coefficients of $w$ at $r^2$ and
$s^2$ must be equal 
for type 2, a contradiction since this forces the $P$-part if the type 1 idempotent to be 0.  
\eop

\begin{coro}\labtt{orbittype1} 
The only idempotents which are a proper summand of some nontrivial  
idempotent are the ones of type  1 and norm $\sixteenth$.  
There are 4 such  and
they come in orthogonal pairs, which are just pair of
idempotents and their conjugates.  
\end{coro}

\begin{coro}  
$Aut(A)$ is a dihedral group of order 8.
\end{coro}  
\pf
The automorphism group preserves and acts faithfully on the
set $J$ of type 1 idempotents of norm $\sixteenth$, the complete set of
idempotents which are proper summands of proper idempotents,  
and furthermore
preserves the partition defined by orthogonality.  
The orthogonal in $A$ of the
nonsingular subspace
$span(J)$ is spanned by $v:=  r^2+s^2 -{ \frac{16+b^2}{4b}}rs$.  
We claim that if an automorphism acts trivially on 
$span(J)$, it acts trivially on $A$.  This is so because 
$rs \in span\{ r^2\times s^2 , r^2, s^2 \}$ and $\{rs\}\cup J$ spans $A$.  

This proves that the automorphism group of $A$ 
embeds in a dihedral group of order 8.   
This embedding is an isomorphism onto since the \lvoap -
group embeds in $Aut(A)$.  
\eop

\begin{prop}\labtt{dih8}  $\avlp$ is just the \lvoap -group, isomorphic to
$Dih_8$.
\end{prop} 
\pf In this case we have $(V_L^+)_2=M(1)^+_2.$ 
Thus any automorphism
of $V_L^+$ preserves $M(1)^+_2.$ 
Now use \ref{lema}.  
\eop

\section{Automorphism group of $V_L^+$ with $rank (L)=2$}

In this section, we assume that the rank of $L$ is equal to 2. If
$L_1=L_2=\emptyset$,  the 
automorphism group of $V_L^+$ was determined 
in Proposition \ref{hrank}. So in this section we assume that
 $L_1$ or $L_2$ is not empty.

\subsection{$L_1=\emptyset$ and $rank (L_2)=2$; $b=0$.   }

Note that $L$ is generated by $L_2.$ 
We will discuss the automorphism group according to the value
$b$ in the Gram matrix $G$ (see \ref{rank02}).

First we assume that $b=0$ in the Gram matrix $G.$ Then $L\cong \sqrt{2}L_{A_1}
\perp  \sqrt{2}L_{A_1}$,   where $L_{A_1}$ is  the
root lattice of type
$A_1.$ Let $L=\ZZ\alpha_1\oplus \ZZ\alpha_2$ with $(\a_i,\a_j)=4\delta_{ij}$
for $i,j=1,2.$ Set 
\begin{eqnarray*}
& &\omega_{1} =\frac{1}{16}\alpha_1(-1)^2
+\frac{1}{4}(e^{\alpha_1}+e^{-\alpha_1}),\nonumber \\
& &\omega_{2}=\frac{1}{16}\alpha_2(-1)^2 -\frac{1}{4}(e^{\alpha_2}+e^{-\alpha_2}).
\end{eqnarray*}
We also use $\a_2$ to define $\omega_3$ and $\omega_4$ in the
same fashion. Then $\omega_i$ for $i=1,2,3,4$ are commutative Virasoro
vectors of central charge $\half$ (see [DMZ] and [DGH]). It is well-known
that $(V_L^+)_2$ is a commutative (nonassociative) algebra under
$u\times v=u_1v$ 
since the degree 1 part is 0 (cf. [FLM]).  
Let  $X$ be the span of
$\omega_i$ for all $i.$

\begin{lem}\labtt{orth} If $u\in (V_L)_2$  is a Virasoro vector of central
charge $1/2$ then $u=\omega_i$ for some $i.$ 
\end{lem}
\pf The space $(V_L)_2$ is 5-dimensional with a basis 
$$\{\omega_1,\omega_2,
\omega_3,\omega_4, \a_1(-1)\a_2(-1)\}.$$ Let $u=\sum_{i=1}^4c_i\omega_i
+x\a_1(-1)\a_2(-1) \in (V_L)_2$ be a Virasoro vector of central charge
$1/2.$ Then $u\times u=2u.$
Note that $\omega_i\times \omega_j=\delta_{i,j}2\omega_i$ for
$i,j=1,2,3,4,$ $\omega_i\times \a_1(-1)\a_2(-1)=\frac{1}{2}\a_1(-1)\a_2(-1)$
and $\a_1(-1)\a_2(-1)\times \a_1(-1)\a_2(-1)=4\a_1(-1)^2+4\a_2(-1)^2.$
So we have a nonlinear system
\begin{eqnarray*}
& & 2c_i=2c_i^2+32x^2, \ i=1,2,3,4\\
& & 2x=\sum_{i=1}^4xc_i.
\end{eqnarray*}
If $x\ne 0$ then $\sum_{i=1}^4c_i=2$ and $2=\sum_{i=1}^4c_i^2+64x^2.$

Since the central charge of $u$ is $1/2$ we have 
$$\frac{1}{4}=u_3u=\sum_{i=1}^4\frac{c_i^2}{4}+16x^2$$
and $1=\sum_{i=1}^4c_i^2+64x^2.$ This is a contradiction. So $x=0.$
This implies that $c_i=0,1$ and $u=\omega_i$ for some $i.$
\eop 

By Lemma \ref{orth}, any automorphism
$\sigma$ of $V_L^+$ induces a permutation of the four 
$\omega_i.$ 

It is known from [FLM] that $(V_L^+)_2$ has a nondegenerate
symmetric bilinear form $(\cdot , \cdot )$ given by
$(u,v)=u_3v$ for $u,v\in (V_L^+)_2.$ The orthogonal 
complement of
$X$ in $(V_L^+)_2$ with respect to the form is spanned by
$\a_1(-1)\a_2(-1).$ Thus $\sigma  \a_1(-1)\a_2(-1)=\lambda \a_1(-1)\a_2(-1)$
for some nonzero constant $\lambda.$ Since
$\a_1(-1)\a_2(-1)\times \a_1(-1)\a_2(-1)= \a_1(-1)^2+\a_2(-1)^2$ which
is a multiple of the Virasoro element $\omega.$ This shows
that $\lambda=\pm 1.$

On the other hand, 
$$V_L^+\cong V_{\sqrt{2}L_{A_1}}^+\otimes V_{\sqrt{2}L_{A_1}}^+
\oplus V_{\sqrt{2}L_{A_1}}^-\otimes V_{\sqrt{2}L_{A_1}}^-.$$
By Corollary 3.3 of [DGH],
$$V_L^+\cong L(\half,0)^{\otimes 4}\oplus L(\half,\half)^{\otimes 4}.$$
So if the restriction of $\sigma$ to $X$ is identity, then
the action of $\sigma$ on $V_{\sqrt{2}L_{A_1}}^+\otimes V_{\sqrt{2}L_{A_1}}^+
$ is trivial and on $V_{\sqrt{2}L_{A_1}}^-\otimes V_{\sqrt{2}L_{A_1}}^-$
is $\pm 1.$ Indeed,  there is  automorphism $\tau$ of $V_L^+$
such that $\tau$ acts trivially on $ V_{\sqrt{2}L_{A_1}}^+\otimes V_{\sqrt{2}L_{A_1}}^+$ and acts as $-1$ on $V_{\sqrt{2}L_{A_1}}^-\otimes V_{\sqrt{2}L_{A_1}}^-$ by the fusion role for $ V_{\sqrt{2}L_{A_1}}^+$ (see [ADL]).  
As $V_{\sqrt{2}L_{A_1}}^+\otimes V_{\sqrt{2}L_{A_1}}^+$ is generated by $\omega_i$ for
$i=1,2,3,4,$ any automorphism preserves
$V_{\sqrt{2}L_{A_1}}^+\otimes V_{\sqrt{2}L_{A_1}}^+$
and its irreducible module 
$V_{\sqrt{2}L_{A_1}}^-\otimes V_{\sqrt{2}L_{A_1}}^-$ (cf. [DM1]).
   As a result,    
$\<\tau\>$ is a normal subgroup of $\avlp$ isomorphic to $\ZZ_2.$ 

Next we show how $Sym_4$ can be realized 
as a subgroup of $Aut (V_L^+)$ by showing that  
any permutation $\sigma\in Sym_4$
gives rise to an automorphism of $V_L^+.$ 
But it is clear that $Sym_4$  acts on $V_L^+$ by permuting the tensor
factors. In order to see that $Sym_4$ acts on $V_L^+$ as automorphisms,
it is enough to show that $\sigma (Y(u,z)v)=Y(\sigma
u,z)\sigma v$ for $\sigma\in Sym_4$ and $u,v\in V_L^+.$ There
are 4 different ways to choose $u,v.$ We only discuss the
case that $u,v\in
 L(\half,\half)^{\otimes 4}$ since the other cases can be 
dealt with in a similar fashion. Let $u=u^1\otimes u^2\otimes u^3\otimes
u^4$ and $v=v^1\otimes v^2\otimes v^3\otimes v^4$ where $u_i,v_i$
are tensor factors in the $i$-th $L(\half,\half).$ Let ${\cal Y}$
be a nonzero intertwining operator of type $\fusion{L(\half,\half)}{L(\half,\half)}{L(\half,0)}.$ Then, up to a constant,
$$Y(u,z)v={\cal Y}(u_1,z)v_1\otimes {\cal Y}(u_2,z)v_2\otimes
{\cal Y}(u_3,z)v_3\otimes {\cal Y}(u_4,z)v_4$$
(see [DMZ]). Since $\sigma$ is a permuation,
it is trivial to verify that $\sigma (Y(u,z)v)=Y(\sigma u,z)\sigma v.$ 

So we
have proved the following:
\begin{prop} If $b=0$ in the Gram matrix $G$ then $L\cong \sqrt{2}L_{A_1}
\times  \sqrt{2}L_{A_1}$ and $\avlp \cong Sym_4\times \ZZ_2.$
\end{prop}

\begin{rem}\labttr{altproofs4x2}  Here is a different proof that 
 $\avlp$ contains a copy of $Sym_4\times \ZZ_2$, using the theory of finite subgroups of Lie groups.  
Our lattice $L$ lies in $M\cong L_{A_1^2}$.   Take $V_M$, which is  a
lattice VOA.  By \cite{DN1} , $V_M$ 
has automorphism group isomorphic to  $PSL(2,\CC) \wr 2$.  In
$PSL(2,\CC)$, there is up to conjugacy a unique four group 
and its normalizer is
isomorphic to
$Sym_4$.  Correspondingly, in $PSL(2,\CC) \wr 2$ 
there is a subgroup isomorphic to 
$Sym_4\wr 2$.  In this, take a subgroup $H$ 
of the form $2^4{:}[Sym_3 \times 2]$. 
Let $t$ be an involution of $H$ which maps to the central involution
of $H/O_2(H) \cong Sym_3\times 2$ and take 
$R:=C_{O_2(H)}(t)\cong 2^2$.    
Take the fixed points $V_M^R$.  We have that  $V_M^R$ is
isomorphic  to our $V_L^+$.  So,   $V_L^+$ gets an action of 
$H/R \cong 2^2{:}[Sym_3 \times 2]\cong Sym_4\times 2$.  
\end{rem}

\subsection{ $L_1=\emptyset$ and $rank (L_2)=2$; $b=2$.  }

Next we assume that $b$ in the Gram matrix 
is 2. Then $L\cong \sqrt{2}L_{A_2}.$  Then
$L=\ZZ\a_1+\ZZ\a_2$ with $(\a_i,\a_i)=4$ and
$(\a_1,\a_2)=2.$ As before we define
$\omega_1,\omega_2$ by using $\a_1,$ $\omega_3,\omega_4$ by
using $\a_2$ and $\omega_5,\omega_6$ by using $\alpha_1+\alpha_2.$ Then
$\omega_i,$ for $i=1,...,6$ form a basis of $(V_L^+)_2.$

\begin{lem}\labtt{allvvb=2} 
 If $u\in (V_L^+)_2$ is a Virasoro vector of central charge
$1/2$, then $u=\omega_i$ for some $i.$
\end{lem}
\pf   {\it First proof. } (There will be a second proof in the next section.)  
 
Let $u=\sum_{i=1}^6 c_i\omega_i$ for some $c_i\in \CC.$
Then $u$ is a Virasoro vector of central charge $1/2$ if and only if
$(u,u)=1/4$ and $u\times u =2u.$   Note that
$$(\omega_i,\omega_i)=1/4, (\omega_{2j-1},\omega_{2j})=0,\ \ 1\leq i\leq 6, j=1,2,3$$
$$(\omega_1,\omega_k)=(\omega_2,\omega_k)=\frac{1}{32}, \ \ k=3,4,5,6.$$
So we have 
\begin{eqnarray*}
(u,u)=\frac{1}{4}\sum_{i=1}^6c_i^2+\frac{1}{16}\sum_{j=1,2}\sum_{j<k\leq 3}
(c_{2j-1}c_{2k-1}+c_{2j-1}c_{2k}+c_{2j}c_{2k-1}+c_{2j}c_{2k})=\frac{1}{4}.
\end{eqnarray*}

In order to compute $u\times u$ we need the following 
multiplication table in $(V_L^+)_2:$ 
$$\omega_{2i-1}\times \omega_{2i}=0, i=1,2,3$$
$$\omega_1\times \omega_3=\frac{1}{4}(\omega_1+\omega_3-\omega_6), \ \ \ 
\omega_2\times \omega_3=\frac{1}{4}(\omega_2+\omega_3-\omega_5)$$
$$\omega_1\times \omega_4=\frac{1}{4}(\omega_1+\omega_4-\omega_5), \ \ \ 
\omega_2\times \omega_4=\frac{1}{4}(\omega_2+\omega_4-\omega_6)$$
$$\omega_1\times \omega_5=\frac{1}{4}(\omega_1+\omega_5-\omega_4), \ \ \ 
\omega_2\times \omega_5=\frac{1}{4}(\omega_2+\omega_5-\omega_3)$$
$$\omega_1\times \omega_6=\frac{1}{4}(\omega_1+\omega_6-\omega_3), \ \ \ 
\omega_2\times \omega_6=\frac{1}{4}(\omega_2+\omega_6-\omega_4)$$
$$\omega_3\times \omega_5=\frac{1}{4}(\omega_3+\omega_5-\omega_2), \ \ \ 
\omega_4\times \omega_5=\frac{1}{4}(\omega_4+\omega_5-\omega_1)$$
$$\omega_3\times \omega_6=\frac{1}{4}(\omega_3+\omega_6-\omega_1), \ \ \ 
\omega_2\times \omega_6=\frac{1}{4}(\omega_2+\omega_6-\omega_2).$$
Then $u\times u=2u$ if and only if 
$$c_1^2+\frac{1}{4}(c_1c_3+c_1c_4+c_1c_5+c_1c_6-c_3c_6-c_4c_5)=c_1$$
$$c_2^2+\frac{1}{4}(c_2c_3+c_2c_4+c_2c_5+c_2c_6-c_3c_5-c_4c_6)=c_2$$
$$c_3^2+\frac{1}{4}(c_1c_3+c_2c_3+c_3c_5+c_3c_6-c_1c_6-c_2c_5)=c_3$$
$$c_4^2+\frac{1}{4}(c_1c_4+c_2c_4+c_4c_5+c_4c_6-c_1c_5-c_2c_6)=c_4$$
$$c_5^2+\frac{1}{4}(c_1c_5+c_2c_5+c_3c_5+c_4c_6-c_1c_4-c_2c_3)=c_5$$
$$c_6^2+\frac{1}{4}(c_1c_6+c_2c_6+c_3c_6+c_4c_6-c_1c_3-c_2c_4)=c_6.$$
There are 
exactly 6 solutions to this linear system: $c_i=1$ and
$c_j=0$ if
$j\ne i$ where $i=1,...,6.$   We thank Harm Derksen for 
obtaining this result with the MacCauley software package.  
This finishes the proof of the
lemma.
\eop  

\begin{prop} If $b=2$ in the Gram matrix, then 
$L\cong \sqrt{2}L_{A_2}$ and $\avlp$ is the LVOA${}^+$-group.  
\end{prop}
\pf First note that 
the Weyl group  acts  on $L$, 
preserving and acting as $Sym_3$ on the set  
$$\{\{\pm \alpha_1\},\{\pm \a_2\},\{\pm(\a_1+\a_2)\}\}.$$

Now let $\sigma\in Aut (V_L^+).$ Set
$X_i=\{\omega_{2i-1},\omega_{2i}\}$ for
$i=1,2,3.$ Then $X_i$ are the only orthogonal pairs in 
$X=X_1\cup X_2\cup X_3.$ Since $\sigma X=X$ we see that $\sigma$ induces 
a permutation on the set $\{X_1,X_2,X_3\}.$ 

The above shows that $\OO (\hat L)$ induces
$Sym_3$ on this 3-set.  We may therefore assume that $\sigma$ preserves
each $X_i.$ In this  case $\sigma$ acts trivially on
$\a_1(-1)^2,\alpha_2(-1)^2, (\a_1+\a_2)(-1)^2.$ That is,
$\sigma$ acts trivially on the subVOA they generate, which
is isomorphic to 
$M(1)^+.$ As a result,
$\s$ is in the LVOA$^+$-group. 
\eop

\subsection{Alternate proof for $b=2$}

The system of equations in the variables $c_i$ which occurred in 
the proof of \ref{allvvb=2} 
can be replaced by an
equivalent system \ref{coeffsw} 
 which looks more symmetric.  The  old system was solved with
software package MacCauley but not with Maple.  The new
system was solved with Maple and gives the same result as
before.

\begin{nota}\labttr{newb=2} 
Let $r$ and $s$ be independent norm 4 elements 
so that $t:=-r-s$ has norm
4. Let $w$ be an idempotent $w=p+q$, where $p=ar^2+bs^2+c
t^2$ and $q=dv_r+ev_s+ev_t$ which satisfies
$(w,w)=\sixteenth$.   
Since $(L,L)\le 2\ZZ$, we may and do assume that the
epsilon-function is identically 1.  
It follows that $v_r\times v_s=v_t$ and similarly for all
permutations of $\{r,s,t\}$.  
\end{nota}  

\begin{lem}\labtt{coeffsw} 
From $w^2=w$, we have equations 
$$ a=16a^2+4ab+4ac-4bc+d^2  \leqno (e1)$$
$$ b=16b^2+4bc+4ba-4ac+e^2  \leqno (e2)$$
$$ c=16c^2+4cq+4cb-4ab+f^2 \leqno (	e3)$$
$$ d=2d(16a+4b+4c)+2ef \leqno (e4)$$
$$ e=2e(4a+16b+4c)+2df \leqno (e5)$$
$$ f=2f(4a+4b+16c)+2de \leqno (e6)$$
and from $(w,w)=\sixteenth$, we get the equation
$$\sixteenth = 32(a^2+b^2+c^2)+16(ab+ac+bc)+2(d^2+e^2+f^2).\leqno (e7)$$
\end{lem}  
\pf  Straightforward from Appendix: Algebraic rules.   
\eop

\begin{prop} \labtt{6sols} 
There are just 6 solutions $(a,b,c,d,e,f)\in \CC^6$ to
the  equations $(e1),\dots, (e7)$.  They are $({1\over 32},0,0,\eighth , 0 ,0  
0) , ({1\over 32},0,0, -\eighth , 0 ,0  
0)$ and ones obtained from these by powers of the permutation $(abc)(def)$.  
\end {prop}  
\pf 
This follows from use of the  {\tt solve} command in the software package Maple.  
\eop 

\begin{rem}\labttr{omite7}   If we omit (e7), there are infinitely many solutions
with
$d=e=f=0$.
The reason is that the Jordan algebra of symmetric degree 2
matrices has infinitely many idempotents.  It seems possible
that the system in Lemma \ref{coeffsw} could be solved by hand.  
\end{rem}

\subsection{ $L_1=\emptyset$ and $rank (L_2)=2$; $b=1$.  }

We now deal with the cases $b=1$ in the Gram matrix.

\begin{prop} If $b=1$ in the Gram matrix, then  $\avlp$ is the LVOA${}^+$
group.  
\end{prop}
\pf By Corollary \ref{orbittype1}, any automorphism of $V_L^+$ 
preserves $M(1)^+_2,$ the result follows from Lemma \ref{lema}
\eop

\section{$\avlp$, for $L_1=\emptyset,$ $rank (L_2)=1$} 

In this case we can assume that $L_2=\{2\alpha_1, -2\alpha_1\}.$ Let
$\alpha_2\in H$ such that $(\alpha_i,\alpha_j)=\delta_{i,j}.$ Then 
$(V_L^+)_2$ is 4-dimensional with basis $v_{2\alpha_1}, \frac{1}{2}\a_1(-1)^2,\frac{1}{2}\a_2(-1)^2, \a_1(-1)\a_2(-1).$

\begin{lem} Any automorphism of $V_L^+$ preserves the subspace $S^2H$
of $(V_L^+)_2$ spanned by $\frac{1}{2}\a_1(-1)^2,\frac{1}{2}\a_2(-1)^2, \a_1(-1)\a_2(-1).$ 
\end{lem}
\pf Since Virasoro vectors of central charge 1 in $S^2H$ span $S^2H,$ 
it is enough to show that any Virasoro vector of central charge 1
lies in $S^2H.$

Let $t=d_1\frac{\a_1^2}{2}+d_2\frac{\a_2^2}{2}+d_3v_{2\alpha_1}+d_4\a_1\a_2$
be a Virasoro vector of central charge 1 with $d_3\ne 0$. 
Then we must have
$t\times t=2t$ and $(t,t)=1/2.$ A straightforward computation shows that
$$t\times t=d_1^2\a_1^2+d_2^2\a_2^2+d_3^2(2\alpha_1)^2+d_4^2(\a_1^2+\a_2^2)
+4d_1d_3v_{2\alpha_1}+2d_1d_4\a_1\a_2+2d_2d_4\a_1\a_2.$$
This gives four equations 
\begin{eqnarray*}
& & d_1=d_1^2+4d_3^2+d_4^2\\
& &d_2=d_2^2+d_4^2\\
& &d_3=2d_1d_3\\
& & d_4=d_1d_4+d_2d_4.
\end{eqnarray*}
The relation $(t,t)=1/2$ gives one more equation:
$$\frac{1}{2}=\frac{1}{2}d_1^2+\frac{1}{2}d_2^2+d_4^2+2d_3^2.$$
Thus 
$$1=d_1+d_2.$$
Since $d_3\ne 0$, $d_1=1/2$ and $d_2=1/2.$ So we have
$$\frac{1}{4}=4d_3^2+d_4^2,  \frac{1}{4}=2d_3^2+d_4^2.$$
This forces $d_3=0,$ a contradiction.
\eop

\begin{prop} In this case,  $\avlp$ is the LVOA${}^+$-group.  
\end{prop}
\pf 
 By Corollary \ref{orbittype1}, any automorphism of $V_L^+$ 
preserves $M(1)^+_2,$ the result follows from Lemma \ref{lema}
\eop 

\section{$\avlp$, for $L_1\ne \emptyset$}

Finally we deal with the case that $L_1\ne \emptyset.$ There are
two cases: rank$(L_1)=2$ or rank$(L_1)=1.$ 
\subsection{$rank (L_1)=2$}
In this case $L=L_{A_1^2}$ or $L=L_{A_2}$ because these are the only rank 2 root lattices possible and each is a maximal even integral lattice in its rational span.

\subsubsection{$L$ has type ${A_1^2}$ } 

If  $L\cong L_{A_1^2}$ then 
$Aut(V_L)\cong PSL(2,\CC)\wr 2$ and 
$$V_L^+ \cong V_{L_{A_1}}^+\otimes V_{L_{A_1}}^+\oplus V_{L_{A_1}}^-\otimes V_{L_{A_1}}^-.$$

Since the connected component of 
the identity in 
$\avl$ contains a lift of
$-1_L$, we may assume that such a 
lift is in a given maximal
torus, so is equal to the 
automorphism 
$e^{\pi i\b(0)/2}$, where $\b$ is a sum of orthogonal
roots.   

It follows that $\vlp \cong V_K$, where $K=2L+\ZZ \b$.
The result \cite{DN1} implies that 
$\avlp \cong  Aut(V_K)$, which is the LVOA group
$\TT_2.Dih_8$.

\subsubsection{ $L$ has type ${A_2}$  } 

Here,  $(V_L^+)_1$ is a 3-dimensional Lie
algebra 
isomorphic to $sl(2,\CC).$ The difficult part in this case is to determine
the vertex operator subalgebra generated by $(V_L^+)_1.$ Let
$L_{A_2}=\ZZ\alpha_1+\ZZ\alpha_2$ such that $(\a_i,\a_i)=2$ and $(\a_1,\a_2)=-1.$
The set of roots in $L$ is $L_1=\{\pm \a_i|i=1,2,3\}$
where
$\alpha_3=\alpha_1+\alpha_2$.  The positive roots are $\{\a_i|i=1,2,3\}.$ The space 
$(V_L^+)_1$ is 3-dimensional with a basis $v_{\a_i}$ 
for  $i=1,2,3$ and  $(V_L^-)_1$ is 5-dimensional with a basis
$\a_1(-1),\a_2(-1), e^{\a_i}-e^{-\a_i}$ 
for $i=1,2,3.$
It is a straightforward to verify that $(v_{\a_i})_{-1}v_{\a_i}$ for $i=1,2,3$ 
and $\a_i(-1)^2$ for $i=1,2,3$  span the same space.   Thus 
$\omega=\frac{1}{4}\a_1(-1)^2+\frac{1}{12}(\a_1(-1)+2\a_2(-1))^2$
lies in the vertex operator algebra generated by $(V_L^+)_1.$

In order to determine the vertex operator
algebra generated by  $(V_L^+)_1$ we need to recall the standard modules
for affine algebra 
$$A_1^{(1)}=\hat{sl}(2,\CC)=sl(2,\CC)\otimes
\CC[t,t^{-1}]\oplus \CC K$$ (cf. [DL]).   
We use the standard basis $\{\a,x_{\a},x_{-\a}\}$ for $sl(2,\CC)$
such that 
$$[\a,x_{\pm\a}]=\pm 2x_{\pm \a}, [x_{\a},x_{-\a}]=\a.$$
We fix an invariant  symmetric nondegenerate bilinear form on
$sl(2,\CC)$  such that  $(\a,\a)=2.$  
The level $k$ standard $A_1^{(1)}$-modules are parametrized by dominant integral
linear weights $\frac{i}{2}\alpha$ for $i=0,...,l$ such that the highest weight
of the $A_1^{(1)}$-module, viewed as a linear form on 
$\CC\a\oplus{\CC}K\subset\hat{sl}(2,\CC),$ is given by
$\frac{i}{2}\alpha$ and the correspondence
$K\mapsto k.$   Let us denote the corresponding standard 
$A_1^{(1)}$-module by $L(k,\frac{i}{2}\a).$ It is well known that $L(k,0)$ is a
simple rational vertex operator algebra and
$L(k,\frac{i}{2}\alpha)$ for $i=0,...,k$ 
is a complete list
of irreducible $L(k,0)$-modules (cf. [DL], [FZ] and [L2]).
Note that 
$$L(k,\frac{i}{2}\a)=\oplus_{n=0}^{\infty}L(k,\frac{i}{2}\a)_{\l_i+n}$$
where $\l_i=\frac{i(i+2)}{4(k+2)}$ and $L(k,\frac{i}{2}\a)_{\l_i+n}$
is the eigenspace of $L(0)$ with eigenvalue $\l_i+n$ (cf. [DL]).  
In fact, the lowest weight space $L(k,\frac{i}{2}\a)_{\l_i}$
of  $L(k,\frac{i}{2}\a)$ is an irreducible 
$sl(2,\CC)$-module of dimension $i+1.$ 

Since $V_{L_{A_2}}$ is a unitary module for affine algebra $A_2^{(1)}$
(cf. [FK]), 
the vertex operator algebra $V$ generated by $(V_L^+)_1$ is isomorphic to
the standard level $k$ $A_1^{(1)}$-module $L(k,0)$
for some nonnegative integer $k.$ Let $\{v_1,v_2,v_3\}$ be an orthonormal
basis of $sl(2,\CC)$ with respect to the standard bilinear form.
Then $\omega'=\frac{1}{2(k+2)}\sum_{i=1}^3v_i(-1)^2\1\in V$ is the
Segal-Sugawara  Virasoro vector. Let 
$$Y(\omega',z)=\sum_{n\in\ZZ}L(n)'z^{-n-2}.$$
Then 
$$[L(n)-L(n)',u_m]=0$$
for $m,n\in \ZZ$ and $u\in V.$ So $L(-2)-L(-2)'$ 
acts as a constant on $V$ as $V$ is a simple vertex operator algebra.
As a result, $L(-2)-L(-2)'=0$ since the left side is 
both a constant and an operator which shifts degree by 2.  
 The creation axiom for VOAs implies that, $\omega'=\omega.$ 
Since the central charge of $\omega$ is 2, the central charge
$\frac{3k}{2(k+2)}$ of $\omega'$ is also 2. This implies
that $k=4$ and $V\cong L(4,0).$  
Now $V_L^+$ is a $L(4,0)$-module and the quotient module
$V_L^+/V$ has minimal weight (as inherited from 
$V_L^+$) greater than 1. On the other hand,
the minimal weight of the irreducible $L(4,0)$-module 
$L(4,\frac{i}{2}\a)$ is  
 $\frac{i(i+2)}{4(4+2)}$ which is less than 2 for $0 \leq i\leq 4$.  
Since every irreducible is one of these, we conclude $V_L^+=V=L(4,0).$ 
Since $V_L^-$ is an irreducible $V_L^+$-module
with minimal weight 1, we immediately see that
$V_L^-=L(4,2\alpha).$  

So we have proved the following:
\begin{prop} If $rank (L_1)=2$ there are two cases.

(1) If $L=L_{A_1^2}$ then $V_L^+$ is again a lattice vertex operator
algebra  $V_K$ where $K$ is generated by $\beta_1,\beta_2$
with $(\beta_i,\beta_i)=4$ and $(\beta_1,\beta_2)=0.$ 
The automorphism
group of $V_L^+$ is the LVOA${}^+$ group which
is isomorphic to the LVOA-group for lattice $K$.

(2) If $L=L_{A_2},$ then $V_L^+$ is isomorphic to the vertex operator
algebra $L(4,0)$ and $\avlp$ is isomorphic to $PSL(2,\CC)$ which
is the automorphism group of $sl(2,\CC).$ 
\end{prop}

\subsection{$rank (L_1)=1$}

\subsubsection{ $L$ rectangular.   } 
We first assume that $L=\ZZ r+\ZZ s$ such that $(r,r)=2,$ 
$(s,s)\in 6+8\ZZ$ and $(r,s)=0.$ Then $V_L=V_{L_{A_1}}\otimes V_{\ZZ s}$ and
$$V_L^+= V_{L_{A_1}}^+\otimes V_{\ZZ s}^+\oplus V_{L_{A_1}}^-\otimes V_{\ZZ s}^-.$$

\begin{lem}\labtt{sometorus} A group of shape   
$(\CC\beta/\ZZ \frac{1}{4}\beta\cdot\ZZ_2)\times \ZZ_2$
acts on $V_L^+$ as automorphisms. 
\end{lem}
\pf
We have already mentioned that $V_{L_{A_1}}^+$ is isomorphic to $V_{\ZZ\beta}$ for $(\b,\b)=8$ and $V_{L_{A_1}}^-$ is isomorphic to $V_{\ZZ\b+\half \b}$
as $V_{L_{A_1}}^+$-modules. We also know from [DN1] that $Aut ( V_{\ZZ \b})$ is
isomorphic to $\CC\b/(\ZZ \frac{1}{8}\b)\cdot \ZZ_2$ where
the generator of $\ZZ_2$ is induced from the $-1$ isometry of the lattice
$\ZZ\b.$ The action of
$\lambda \beta\in \CC\b$ is given by the operator $e^{2\pi i \lambda\b(0)}.$
Note that $\CC\b$ acts on  $V_{\ZZ\b +\half \b}$
in the same way. But the kernel of the action of $\CC\b$ on 
$V_{\ZZ\b +\half \b}$ is $\ZZ \frac{1}{4}\beta$ instead
of $\ZZ \frac{1}{8}\beta.$ As a result, the torus  
$\CC\beta/\ZZ \frac{1}{4}\beta$ acts on both  $V_{\ZZ\b}$ and $
V_{\ZZ\b +\half \b}.$ 
By [DG],
$Aut (V_{\ZZ s}^+)$ is isomorphic to 
$\frac{1}{2}{\ZZ s}/\ZZ s\cong\ZZ_2$ which also acts on 
$V_{\ZZ s}^-.$ So the group $(\CC\beta/\ZZ \frac{1}{4}\beta\cdot\ZZ_2)\times \ZZ_2$
acts on $V_L^+$ as automorphisms. 
\eop 

In order to determine $Aut(V_L^+)$ in this case we need to recall
the notion of commutant from [FZ].

\begin{de}\labttr{commutant}  Let $V\!=\!(V,Y,\1,\omega)$ be
a vertex operator algebra and $U\!=\!(U,Y,\1,\omega')$ be vertex
operator subalgebra with a different Virasoro vector $\omega'.$ The
{\it commutant}  $U^c$ of $U$ in $V$ is defined by
$$U^c :=\{v\in V|u_nv=0, u\in U,n\geq 0\}.$$
\end{de}

\begin{rem}\labttr{vacuumlike} 
The above space   $U^c$ is  the space of vacuum-like vectors for $U$
(see [L1]).
\end{rem} 

\begin{lem}\labtt{llle} Let $V$ be a vertex operator algebra and $U^i=
(U^i,Y,1,\omega^i)$ are simple vertex operator subalgebras 
of $V$ with Virasoro vector $\omega^i$ for $i=1,2$ such that
$\omega=\omega^1+\omega^2.$ We assume that 
$V$ has a decomposition 
$$V \cong \oplus_{i=0}^pP^i\otimes Q^i$$
as $U^1\otimes U^2$-module such that $P^0\cong U^1,$ $Q^0\cong U^2$,  the 
$P^i$ are inequivalent $U^1$-modules and the $Q^i$ are
inequivalent $U^2$-modules. Then $(U^1)^c=U^2$ and $(U^2)^c=U^1.$
\end{lem}
\pf It is enough to prove that $(U^2)^c\subset U^1.$ Let $v\in (U^2)^c.$
Then $v$ is a vacuum-like vector for $U^2.$ Then the $U^2$-submodule
generated by $v$ is isomorphic to $U^2$ (see [L1]). Since 
$V$ is a completely reducible $U^2$-module and any
$U^2$-submodule isomorphic to $U^2$ is contained in
$U^1\otimes U^2.$ In particular,
$v\in U^1\otimes U^2.$ This forces $v\in U^1.$
\eop

\begin{prop}\labtt{p1} The  group $Aut(V_L^+)$ 
is isomorphic to $((\CC\beta/ \ZZ \frac{1}{4}\beta)\cdot \ZZ_2)\times \ZZ_2.$  This can be interpreted as an action of $\NN (\widehat {\ZZ \b}) \times \ZZ_2$, where $(\b ,\b )=8$.
\end{prop}
\pf We have already shown \ref{sometorus} that
the group $((\CC\beta/\ZZ \frac{1}{4}\beta)\cdot \ZZ_2)\times \ZZ_2$ acts
on $V_L^+$ as automorphisms. 

Let $\sigma$ be an automorphism of $V_L^+.$ Then 
$\sigma \beta(-1)=\l \beta(-1)$ for some nonzero $\l\in \CC$ as $(V_L^+)_1$
is spanned by $\b(-1).$ This implies that $\sigma \beta(n)\sigma^{-1}=
\lambda \beta(n)$ for $n\in\ZZ.$ 
Since $V_{\ZZ s}^+$ is precisely the subspace
of $V_L^+$ consisting of vectors killed by $\beta(n)$ for $n\geq 0$, 
we see that $\sigma   V_{\ZZ s}^+\subset V_{\ZZ s}^+.$ 
Thus $\sigma|_{V_{\ZZ s}^+}$ is an automorphism of $V_{\ZZ s}^+.$
On the other hand, $V_{L_{A_1}}^+$  is the commutant of $V_{\ZZ s}^+$ in $V_L^+$by Lemma \ref{llle}.

The above show that  $\sigma$ induces an automorphism of the tensor factor  $V_{L_{A_1}}^+.$
The restriction of $\sigma$ to $ V_{L_{A_1}}^+\otimes V_{\ZZ s}^+$
is a product  $\sigma_1\otimes \sigma_2$ for
some $\sigma_1\in Aut (V_{L_{A_1}}^+)$ and $\sigma_2\in Aut (V_{\ZZ s}^+).$
Multiplying $\sigma$ by $\sigma_2$ we can assume that $\sigma=1$ on
$V_{\ZZ s}^+.$ As we have already mentioned,  $Aut (V_{L_{A_1}}^+)$
is isomorphic to $(\CC\beta/\ZZ \frac{1}{8}\beta){\cdot} \ZZ_2.$ 
Since $(\CC\beta/\ZZ \frac{1}{8}\beta)$ acts trivially on $\beta(-1)$
and the outer factor $\ZZ_2$ is represented in  $Aut( V_L^+)$ by  
action of $\pm 1$ on $\beta(-1)$.   
As a result $\sigma \beta(-1)=\pm \beta(-1).$ 
Now multiplying $\sigma$ by an outer element of $(\CC\beta/\ZZ \frac{1}{4}\beta)\cdot \ZZ_2$,  we can assume that $\sigma \b(-1)=\beta(-1).$

Set $W_{n\beta}=M(1)\otimes e^{n\beta}\otimes  V_{\ZZ s}^+$ and 
$W_{n\beta+\beta/2}=M(1)\otimes e^{n\beta+\beta/2}\otimes 
V_{\ZZ s}^-$ for $n\in\ZZ$ where $M(1)=\CC[\beta(-n)|n>0].$ 
 Then
$V_L^+=\oplus_{n\in\frac{1}{2}\ZZ}W_{n\beta}$ and $u_mv\in
W_{\mu+\nu}$ for $u\in W_{\mu}$ and $v\in W_{\nu},$ and
$n\in \ZZ.$ Note that $W_{\mu}$ is the eigenspace of
$\beta(0)$ with eigenvalue $(\beta,\mu).$  Since $\sigma
\beta(-1)=\beta(-1)$ we see that
$\sigma$ acts on each $W_{\mu}$ as a constant $\l_{\mu}$ and
$\l_{\mu}\l_{\nu}=\l_{\mu+\nu}.$ As a result, $\sigma=
e^{2\pi i \gamma(0)}$ for some $\gamma\in 
\CC \beta.$ That is, $\sigma$ lies in $\CC\beta/\ZZ \frac{1}{4}\beta.$ 
This completes the proof.
\eop

\subsubsection{ $L$ not rectangular }  
Next we assume that $L\ne \ZZ r\perp  \ZZ s.$ Then $L=\ZZ
r\oplus \ZZ\frac{1}{2}(s+ t)$ where $(s,s)\in 6+8\ZZ$ and
$(s,s)\geq 14$  (see \ref{rank1}).  
 Let $K=\ZZ r\oplus \ZZ
s.$ Then
$L=K\cup (K+\frac{1}{2}(r+s))$ and 
$V_L=V_{\ZZ r}\otimes V_{\ZZ s}\oplus V_{(\ZZ+\frac{1}{2})r}
\otimes
V_{(\ZZ+\frac{1}{2})s}.$ Thus
$$V_L^+= V_{\ZZ r}^+\otimes V_{\ZZ s}^+ \oplus  V_{\ZZ r}^-\otimes V_{\ZZ s}^-
\oplus V_{(\ZZ+\frac{1}{2})r}^+\otimes
V_{(\ZZ+\frac{1}{2})s}^+ \oplus V_{(\ZZ+\frac{1}{2})r}^-\otimes
V_{(\ZZ+\frac{1}{2})s}^-$$ and
$$V_L^+=V_K^+\oplus V_{K+\frac{1}{2}(s+t)}^+.$$
As before,  we note that $V_{\ZZ r}^+$ is
isomorphic to $V_{\ZZ \beta}$ with $(\b,\b)=8.$ 

\begin{prop}\labtt{nonrectangular}  Assume that $rank(L_1)=1,$ $L\ne \ZZ r+\ZZ s$, $r, s$ as above.   Then 
$\avlp \cong (\CC\beta/\frac{1}{2}\ZZ\b)\cdot \ZZ_2$,  where $(\b,\b)=8$.  
The action is trivial on the subVOA $V_{\ZZ s}^+$ and leaves $V_{\ZZ r}^+$ invariant.  
A generator
of the quotient $\ZZ_2$  comes from the $-1$ isometry of $\frac{1}{4}\ZZ\beta$ and
$\alpha\in \CC\beta$ acts as $e^{2\pi i\alpha(0)}.$
\end{prop}

\pf Note that $V_K^+$ is a subalgebra of $V_L^+$ and $V_{K+\frac{1}{2}(r+s)}^+$
is an irreducible $V_K^+$-module. By Proposition \ref{p1},
$$Aut (V_K^+)=((\CC\beta/\frac{1}{4}\ZZ\b){\cdot} \ZZ_2)\times \ZZ_2.$$
As we have already mentioned that $V_{\ZZ r}^+$ is isomorphic to
$V_{\ZZ\b}$ with $(\b,\b)=8$ and $V_{\ZZ r}^-$ is isomorphic to 
$V_{\ZZ\beta+\half \beta}$ as $V_{\ZZ\b}$-module. It is easy
to see that $V_{(\ZZ+\frac{1}{2})\beta}^{\pm}$ is isomorphic
to $V_{(\ZZ\pm\frac{1}{4})\beta}$ as $V_{\ZZ\b}$-module. So the action
of $\CC\beta/\ZZ \frac{1}{4} \b$ on $V_K^+$ cannot be extended to an action
of $V_L^+.$ But the torus   $\CC\beta/\frac{1}{2}\ZZ\b$ does acts on
$V_L^+.$ As a result, $\NN (\widehat {\ZZ \half \b} ) \cong (\CC\beta/\frac{1}{2}\ZZ\b ) \cdot \ZZ_2$
is a subgroup of $\avlp.$

The same argument used in the proof of Proposition 
\ref{p1} shows that that any automorphism $\sigma$ of $V_L^+$ preserves
$V_{\ZZ r}^+\otimes V_{\ZZ s}^+$.    Since 
$V_{\ZZ r}^+\otimes V_{\ZZ s}^+,$ $V_{\ZZ r}^-\otimes V_{\ZZ s}^-$,
$V_{(\ZZ+\frac{1}{2})r}^+\otimes
V_{(\ZZ+\frac{1}{2})s}^+,$ $V_{(\ZZ+\frac{1}{2})r}^-\otimes
V_{(\ZZ+\frac{1}{2})s}^-$ are inequivalent irreducible
$V_{\ZZ r}^+\otimes V_{\ZZ s}^+$-modules (see [DM1] and [DLM]),
we see that $\sigma$ preserves 
$$V_K^+= V_{\ZZ r}^+\otimes V_{\ZZ s}^+ \oplus  V_{\ZZ r}^-\otimes V_{\ZZ s}^-.$$

Since $\CC\beta/\frac{1}{4}\ZZ\b\cdot \ZZ_2$ is a quotient group of 
$\CC\beta/\frac{1}{2}\ZZ\b\cdot \ZZ_2,$ we can multiply $\sigma$ by an element of
$\CC\beta/\frac{1}{2}\ZZ\b\cdot
\ZZ_2$ and assume that $\sigma$ acts trivially on the first tensor factor of
$V_K^+.$ If $\sigma$ is the identity on $V_K^+$, then $\sigma$ is either 1 or $-1$
on $V_{K+\frac{1}{2}(r+s)}^+.$ If $\sigma$ is $-1$ on  $V_{K+\frac{1}{2}(r+s)}^+$ then $\sigma=e^{\pi i\frac{1}{2}\beta(0)}$ is an element 
of  $\CC\beta/\frac{1}{2}\ZZ\b\cdot
\ZZ_2.$

If $\sigma$ is not identity on $V_K^+$ then we must have 
$\sigma=e^{\pi i\frac{1}{(s,s)}s(0)}$ on $V_K^+.$   
We will get a
contradiction in this case. Notice that the lowest weight space of
$V_{(\ZZ+\frac{1}{2})r}^+\otimes V_{(\ZZ+\frac{1}{2})s}^+$ is 1-dimensional and
spanned by
$u=(e^{r/2}+e^{-r/2})\otimes (e^{s/2}+e^{-s/2}).$ Since $\sigma$ preserves
 $V_{(\ZZ+\frac{1}{2})r}^+\otimes
V_{(\ZZ+\frac{1}{2})s}^+$, it must map $u$ to $\lambda u$ for some nonzero
constant $\l.$ Note that $u_{\fourth (r+s,r+s)-1}u=4.$  This forces
$\lambda=\pm 1.$ On the other hand, 
$$u_{-\fourth (r+s,r+s)-1}u=(e^r+e^{-r})\otimes (e^s+e^{-s})+\cdots$$
has nontrivial projection to the  $- 1$ eigenspace  of $\sigma$ in $V_K^+.$  
This forces $\lambda=\pm i, $ a contradiction.  
\eop

\section{Appendix: Algebraic rules}

For the symmetric matrices of degree $n$, there is a widely
used basis, Jordan product and inner product, which we
review here.  (This section is taken almost verbatim from [G4]).  

\begin{prop} 
$H$ is a vector space of finite dimension $n$ with
nondegenerate symmetric bilinear form $(\cdot  , \cdot )$.  

$r, s, \dots $ stand for elements of $H$ and $rs$ stands for
the symmetric tensor $r \otimes s + s \otimes r$.  

$rs \times pq = (r,p)sq+(r,q)sp+(s,p)rq+(s,q)rp$.

$(rs,pq)=(r,p)(s,q)+(r,q)(s,p)$

$rs \times v_t= (r,t)(s,t) v_t$.  
\end{prop}  


\begin{de}\label{sbfv2} \rm 
 {\sl The Symmetric Bilinear Form.}  
Source: [FLM], p.217.  This form is associative with respect to the product
(Section 3).  We write  $H$ for  $H_1$.  The set of all $g^2$ and
$x_\alpha^+$ spans $V_2$.   $$ \langle g^2, h^2 \rangle= 2 \langle g,
h \rangle ^2, \leqno (2.2.1)$$

whence 

$$\langle pq, rs \rangle =   \langle p,r \rangle  \langle q,s \rangle + 
 \langle p,s \rangle  \langle q,r \rangle, \text { for } p,q,r,s \in H. \leqno
(2.2.2)$$

$$ \langle x_\alpha^+, x_\beta^+ \rangle = \begin{cases} 2 & \alpha =
\pm
\beta
\cr  0 & \ else \cr  \end{cases}  
 \leqno (2.2.3)
$$

$$\langle g^2, x_\beta^+ \rangle = 0. \leqno (2.2.4)$$

\end{de}

\begin{de}\label{vireltidentity}\rm  
In addition, we have  the distinguished Virasoro
element $\omega$ and  identity $\II := {1 \over 2} \omega$ on $V_2$ (see
Section 3).    If $h_i$ is a basis for $H$ and $h_i^*$ the dual basis, then
$\omega = {1 \over 2} \sum_i h_ih_i^*$.  
\end{de}

\begin{rem}\label{someips}\rm

$$\langle g^2, \omega \rangle = \langle g,g \rangle \leqno (2.4.1)$$

$$\langle g^2, \II  \rangle = {1 \over 2}    \langle g,g \rangle
\leqno (2.4.2)$$

$$\langle   \II ,  \II \rangle = dim(H)/8 
\leqno (2.4.3)$$

$$\langle \omega, \omega\rangle = dim(H)/2  \leqno (2.4.4)$$

If $\{x_i~|~i=1,\dots \ell\}$ is an ON basis, 

$$\II = {1 \over 4} \sum_{i=0}^\ell x_i^2 \leqno (2.4.5)$$

$$\omega = {1 \over 2}  \sum_{i=0}^\ell x_i^2. \leqno (2.4.6)$$

\end{rem}

\begin{de}\label{productv2} \rm 
   {\sl The product on $V_2^F$}  comes from the vertex
operations.  We give it on standard basis vectors, namely $xy \in
S^2H_1$, for $x, y \in H_1$ and $v_\lambda := e^\lambda +
e^{-\lambda}$, for $\lambda \in L_2 $.  (This is the same as $x_\l^+$, used in [FLM].)  Note that (3.1.1) give the
Jordan algebra structure on $S^2H_1$, identified with the space of
symmetric  $8 \times 8$ matrices, and with $\la x,y \ra = {1 \over
8}tr(xy)$.   The function $\varepsilon$ below is a standard part of
notation for
 lattice VOAs.  

$$ x^2 \times y^2 = 4\la x,y \ra xy,  \hskip 1cm  pq \times
y^2=2\la p,y \ra qy+2\la q,y \ra py,  \leqno  (3.1.1)   $$  
$$pq \times rs=\la p,r \ra qs+\la p,s \ra qr+\la q,r \ra ps+\la q,s \ra pr;     
 $$

$$x^2 \times v_\lambda = \la x,\lambda \ra ^2
v_\lambda, 
\hskip 1cm xy \times v_\lambda = \la x, \lambda \ra \la y, \lambda \ra
v_\lambda   \leqno   (3.1.2)$$ 
$$v_\lambda \times v_\mu = 
\begin{cases} 0 &  \la \lambda,\mu \ra  \in \{0,\pm 1,\pm 3\}; \cr
\varepsilon \la \lambda,
\mu \ra  v_{\lambda + \mu} &  \la \lambda,\mu \ra  =-2; \cr \lambda^2 &
\lambda=\mu.\end{cases} 
  \leqno (3.1.3) $$
\end{de} 


Some consequences are these: 

\begin{coro}\labtt{formulaforidentity}  
If $x_1, \dots $ is a basis and $y_1, \dots $ is the dual
basis, then $\II  := \fourth \sum_{i=1}^n x_iy_i$ is the
identity of the algebra $S^2H$.

$(\II, \II) = \frac{n}{8}$.  
\end{coro}


\begin{thebibliography}{10000}
\bibitem[AD]{AD} T. Abe and C. Dong, Classification of irreducible modules
for vertex operator algebra $V_L^+:$ general case, {\em J. Algebra}
{\bf  273} (2004), 657--685.
\bibitem[ADL]{ADL} T. Abe, C. Dong and H. Li, Fusion rules for the vertex 
operator algebras $M(1)^+$ and $V_L^+,$ {\em Comm. Math. Phys.,} to appear.
 
\bibitem[B]{B} R. E. Borcherds, Vertex algebras, Kac-Moody algebras, and the Monster,
{\it Proc. Natl. Acad. Sci. USA} {\bf 83} (1986), 3068-3071.

\bibitem[D]{D} C. Dong, Vertex algebras associated with even lattices, 
{\em J. Algebra} {\bf 161} (1993), 245-265.

\bibitem[DG1]{DG1} C. Dong and R.L. Griess Jr.,
 Rank one lattice type vertex operator algebras and their
automorphism groups, {\em J. Algebra} {\bf 208} (1998), 262-275.

\bibitem[DG2]{DG2} C. Dong and R.L. Griess Jr.,
Automorphism groups of  finitely generated vertex operator 
algebras, {\em Michigan Math. J.} {\bf 50} (2002), 227--239.

\bibitem[DGR]{DGR} C. Dong and R.L. Griess Jr. and A. J. E. Ryba,
 Rank one lattice type vertex operator algebras and their
automorphism groups,  {\em J. Algebra} {\bf  217} (1999), 701--710. 

\bibitem[DGH]{DGH} 
C.~Dong, R.~Griess Jr. and G.~Hoehn,
Framed vertex operator algebras, codes and the moonshine module,
{\em Comm. Math. Phys.} {\bf 193} (1998), 407--448.

\bibitem[DL]{DL}
C.~Dong and J.~Lepowsky, \textit{Generalized vertex algebras and relative vertex operators}, 
Progress in Math., Vol.{\bf112}, Birkh{\"a}user, Boston,
1993.


\bibitem[DLM]{DLM} C. Dong, H. Li and G. Mason,
Compact automorphism groups of vertex operator algebras,  {\em
International Math. Research Notices} {\bf 18} (1996), 913--921.

\bibitem[DM1]{DM1} C. Dong and G. Mason, 
On quantum Galois theory,  {\em Duke Math. J.} {\bf 86}
(1997), 305--321.

\bibitem[DM2]{DM2} C. Dong and G. Mason, Vertex operator algebras and their automorphism groups,
In: Proceedings of International Conference
on Representation Theory (Shanghai, 1998), China Higher Education
Press and Springer-Verlag, Beijing, 2000, 145-166.

\bibitem[DN1]{DN1} C. Dong and K. Nagatomo, Automorphism groups
of lattice vertex algebras,  {\em Contemp. Math.} {\bf 248} (1999), 117-133.

\bibitem[DN2]{DN2} C. Dong and K. Nagatomo, Classification of irreducible modules for the vertex operator algebra $M(1)^+$ II. Higher Rank, \textit{J.~Algebra}, {\bf240} (2001), 389--325.

\bibitem[FLM]{FLM}
I. Frenkel, J. Lepowsky and A. Meurman, {\it Vertex Operator Algebras
and the Monster}, Pure and Appl. Math., {\bf Vol. 134}, Academic Press,
Boston, 1988.

\bibitem[FZ]{FZ}
I.~Frenkel and Y.~Zhu, Vertex operator algebras associated to representations of
affine and Virasoro algebras, {\it Duke~Math.~J.} {\bf66} (1992), 123--168.

\bibitem[G1]{G1} R. Griess, Jr.,  GNAVOA, I. Studies in groups, nonassociative algebras and vertex operator algebras. Vertex operator algebras in 
mathematics and physics (Toronto, ON, 2000), 71--88, Fields Inst. Commun., {\bf39}, AMS, 2003.

\bibitem[G2]{G2} Robert L. Griess, Jr., Automorphisms of vertex operator algebras, a survey, 
Proceedings of the Raleigh Conference on affine algebras, quantum affine
algebras and related topics, 21-24 May, 1998.  
Contemporary Mathematics, volume 248,  American Mathematical Society,  
1999.  

\bibitem[G3]{G3}  Robert L. Griess, Jr.,   Positive definite
lattices of rank at most 8,  Journal of Number Theory, 103
(2003), 77-84. 

\bibitem[G4]{G4} R. Griess, A vertex operator algebra related to
$E_8$ wth automorphism group $O^+(10,2)$, article inThe Monster
and Lie Algebras, ed. J. Ferrar and K. Harada, de Gruyter, 1998.  

\bibitem [GH]{GH} R.~Griess Jr. and G.~Hoehn, Frame stabilizers for
 the lattice vertex operator algebra of  type $E_8$, 
with Gerald H\"ohn, {\em J. Reine. Angew. Math.} {\bf 561} (2003), 1-37. 

\bibitem[J]{J} 
Nathan Jacobson,  
"Isomorphism of Jordan Rings", Amer. J. Math. 70,
(1948), 317-326.
 

\bibitem[L1]{L1} H.~Li, Symmetric invariant bilinear forms on vertex operator algebras,
{\it  J. Pure Appl. Algebra} {\bf  96} (1994), 279--297. 

\bibitem[L2]{L2} H.~Li, Local systems of vertex operators, vertex
superalgebras and modules, \textit{J.~Pure~Appl.~Alg}, {\bf109}
(1996), 143--195.

\bibitem[S]{S} H. Shimakura, The automorphism group of the vertex operator algebra $V_L^+$ for an even lattice $L$ without roots, {\em J. Algebra} {\bf 280} (2004), 29-57.


\end{thebibliography}
\end{document}